\documentclass[preprint,11pt]{elsarticle}

\usepackage{amssymb}
 \usepackage{amsmath}
\usepackage[colorlinks]{hyperref}
\usepackage{doi}
 \newdefinition{remark}{Remark}
 \newtheorem{proposition}{Proposition}
 \newtheorem{theorem}{Theorem}
 \newproof{proof}{Proof}
 \usepackage{algorithm}
 \usepackage{algpseudocode}

\journal{}

\begin{document}

\begin{frontmatter}

\title{Numerical Approximation of Optimal Convex and Rotationally Symmetric Shapes for an Eigenvalue Problem arising in Optimal Insulation}

\author[HK]{Hedwig Keller\corref{cor1}}
\ead{hedwig.keller@mathematik.uni-freiburg.de}
\author[SB]{Sören Bartels}
\ead{bartels@mathematik.uni-freiburg.de}
\author[GW]{Gerd Wachsmuth}
\ead{gerd.wachsmuth@b-tu.de}

\affiliation[HK]{organization={Albert-Ludwigs-University Freiburg,Department of Applied Mathematics},
 addressline={Hermann-Herder-Str. 10},
             city={Freiburg},
             postcode={79104},
             country={Germany}}

\affiliation[SB]{organization={Albert-Ludwigs-University Freiburg, Department of Applied Mathematics},
 addressline={Hermann-Herder-Str. 10},
             city={Freiburg},
             postcode={79104},
             country={Germany}}

\cortext[cor1]{Corresponding author}

\affiliation[GW]{organization={BTU Cottbus-Senftenberg,  Optimale Steuerung},
            ardressline =  {Postfach 10 13 44},
            city = {Cottbus},
            postcode = {03013},
            country = {Germany}}

\begin{abstract}

We are interested in the optimization of convex domains under a PDE constraint.
Due to the difficulties of approximating convex domains in $\mathbb{R}^3$, the
restriction to rotationally symmetric domains is used to reduce shape
optimization problems to a two-dimensional setting. For the optimization of an
eigenvalue arising in a problem of optimal insulation, the existence of an
optimal domain is proven. An
algorithm is proposed that can be applied to general shape optimization
problems under the geometric constraints of convexity and rotational symmetry. The approximated optimal domains
for the eigenvalue problem in optimal insulation are discussed.

\end{abstract}

\begin{keyword}
shape optimization \sep optimal insulation  \sep convexity \sep rotational symmetry \sep PDE constraints   \sep iterative solution

\MSC[2020] 49Q10 \sep 49M41 \sep 65N25

\end{keyword}

\end{frontmatter}

\section{Introduction}

Solvability of shape optimization problems relies, among other factors, on
strong constraints on the geometry of the admissible domains.
Since we minimize over shapes, no topology is readily available. The
restriction to classes of convex domains appears attractive, since the
compactness results available for convex domains let us avoid more general
topological frameworks. For corresponding analytical details we refer to
\cite{variationalmethods}, \cite{veodf}, \cite{vangoethem},
\cite{yang09}, \cite{BD03} and \cite{BG}.
Therefore,
we restrict the shape optimization to
 open, convex and bounded domains. 
 
 However, numerical approximation of convex domains is difficult in higher
 dimensions. Indeed, for conformal P1 finite elements we can not guarantee that
 a convex function can be approximated consistently (c.f. \cite{chone}), and with simple
 examples we can show, that the nodal interpolant of a convex function is not
 necessarily convex itself, for such an example see \cite[Figure 2.1]{aguilera}.
 To approximate convex functions,
 we need for example higher order conforming finite elements (c.f.
 \cite{wachsmuth}), a weaker definition for convexity tailored to finite
 elements (c.f. \cite{aguilera}), a geometric approach as in
 \cite{convexhull} or spherical harmonic decomposition (c.f.
 \cite{antunes}).
 Since the approximation
 of convex domains in $\mathbb{R}^3$
 has certain similarities to the approximation
 of convex functions in $\mathbb{R}^2$,
 we expect related difficulties.
Therefore, we restrict our domains to a class of rotationally symmetric
domains, which allows us to reduce the problem to a two-dimensional setting,
for which the boundary is a convex curve. The dimensional reduction also allows for a higher resolution in the numerical approximation.

We are interested in the optimization under a PDE constraint, in particular in optimizing an eigenvalue occurring in a problem of optimal insulation. For more details in PDE constraint optimization we refer to \cite{hinze}.

A heat conducting body is to be coated by an insulating
material in such a way to get the best insulating properties. This translates
to the non-linear eigenvalue problem
\[
\lambda_m = \min \left\lbrace J_m(u) := \int_{\Omega} \vert \nabla u\vert^2 \, \text{d}x + \frac{1}{m}\left( \int_{\partial \Omega} \vert u \vert \, \text{d}s \right)^2 \ : \ \int_{\Omega} \vert u \vert ^2 \, \text{d}x = 1  \right\rbrace.
\]
From \cite{BBN17} we expect that in general the distribution of insulating
material is asymmetric and that the ball is not optimal, in contrast to what we
might expect from isoparametric inequalities for eigenvalues of the Laplacian. 

The numerical framework for the approximation of the eigenvalue from
\cite{BB19} confirmed the expected asymmetry in two dimensions.  Our goal
is to perform the shape optimization for convex, rotationally symmetric domains
in $\mathbb{R}^3$. The numerical experiments in Section  \ref{chap:c_num_ex} confirm, that the constraint to rotational symmetric domains and eigenfunctions still allows for a break in symmetry.

We focus on the existence of an optimal domain and the meaningful numerical approximations provided by the proposed algorithm. We will discuss the stability of the numerical scheme shortly, but a detailed examination lies beyond the scope of this work. In the proof of existence the geometric constraints, especially the convexity, play key roles. 

First, in Section \ref{sec:rot_sym} we describe the dimensional reduction
obtained from the rotational symmetry. Then we consider the shape
optimization for the eigenvalue problem arising in the problem of optimal
insulation. We prove existence of an optimal domain in Section \ref{chap:ex_theory} and derive the two-dimensional problem and its numerical approximation
and comment on the stability of the numerical scheme in Section \ref{chap:dis_red_p}. In
Section \ref{chap:num_approx} we establish a framework for the numerical approximation of optimal convex domains described in \cite{BW18} but
adjusted for rotational symmetry, which can be applied to different shape
optimization problems as well. The numerical experiments are evaluated in
Section \ref{chap:c_num_ex}.

\section{Rotationally Symmetric Domains and Dimensional Reduction}
\label{sec:rot_sym}

We consider a shape optimization problem that, for a given open and
bounded domain $\widehat{Q}$, density function $j$, volume $M$ and state
equation $\widetilde{\varphi}$, seeks a domain $\Omega$  which solves
\begin{align*} \label{problem}
\text{ Minimize } & \ \int_{\Omega} j(x,u(x),\nabla u(x))\, \text{d}x \\
\text{w.r.t }&\ \Omega \subset \widehat{Q} \subset \mathbb{R}^3  \text{ open, convex and rotationally symmetric} \\ \text{and } & \int_\Omega \, \text{d}x = M  \tag{\textbf{P}}\\
\text{s.t. }&\ u \in H^1(\Omega) \text{ solves a certain state equation } \widetilde{\varphi}(u) = 0.
\end{align*}
Here,
the rotational symmetry is to be understood w.r.t.\ the $x_3$-axis.
We assume that $\widehat Q$ and $j$ are rotationally symmetric as well.  Furthermore, we assume, that the solution $u \in H^1(\Omega)$ of the state equation is rotationally symmetric, based on analytic properties or results of numerical experiments of the problems under consideration.
For example, in the eigenvalue problem considered in Section \ref{chap:ex_theory}, previous experiments suggests that the eigenfunctions of the ball are rotationally symmetric, c.f. \cite{BB19}.

We use the rotational symmetry to reduce the problem to a
two-dimensional setting. For this, we first use a transformation to cylindrical
coordinates and then neglect the angle due to it being
constant because of the rotational symmetry.
For one half of the cross section $Q \subset \mathbb{R}^2$  of $\widehat{Q}$ we define $\Phi :Q \times [0,2\pi) \rightarrow \widehat{Q}$ as the transformation
from cylindrical coordinates  to Cartesian coordinates $\Phi(r, z, \phi) = (r
\cos \phi, r \sin \phi, z) = (x_1, x_2, x_3)$.

From now on, we only consider rotationally symmetric functions, i.e.
functions from the space 
\begin{equation*}
H^1_{\textup{sym}}(\widehat{Q}) := \{ \widehat{u} \in H^1(\widehat{Q}) :   \partial_\phi \widehat{u}= 0\}.
\end{equation*}
Since the set of rotationally symmetric functions is closed under $L^1$-convergence,
$H^1_{\textup{sym}}(\widehat{Q})$ with the $H^1$-norm is also a Hilbert space.
For a function $\widehat{u} \in H^1_{\textup{sym}}(\widehat{Q})$ we then associate the dimensionally reduced function as
\begin{equation*}
u(r,z) = \frac{1}{2\pi}\int_0^{2\pi} \widehat{u} \circ \Phi(r,\phi,z) \, \text{d}\phi  \quad \text{for }(r,z) \in Q.
\end{equation*}
We show that $u$ is also weakly differentiable with regard to the variables $(r,z)$.  For a test function $\varphi \in C^\infty_c(Q,\mathbb{R}^2)$ we have 

\begin{align*}
-\int_Q  u \text{ div}\varphi\text{d}x &= -\int_Q \frac{1}{2\pi} \int_0^{2\pi} (\widehat{u} \circ \Phi)(r,\phi,z)  \text{ div}\varphi \, \text{d}\phi \text{d}(r,z) \\ &= \int_Q \frac{1}{2\pi}\int_0^{2\pi} D\Phi_{r,z}(r,\phi,z)^\top \nabla\widehat{u}(\Phi(r,\phi,z)) \cdot \varphi \, \text{d}\phi \, \text{d}(r,z) \\ &= \int_Q  \nabla u  \cdot \varphi \, \text{d}(r,z)
\end{align*}
with $D\Phi_{r,z}(r,\phi,z)$ the Jacobi matrix with respect to  $r$ and $z$ and the weak derivative 
\begin{equation*}
 \nabla u = \frac{1}{2\pi} \int_0^{2\pi}  D\Phi_{r,z}(r,\phi,z)^\top \nabla\widehat{u}(\Phi(r,\phi,z)) \, \text{d}\phi.
\end{equation*}

We take a closer look at the relation between the functions $\widehat{u}$
of  $H^1_{\textup{sym}}(\widehat{Q})$ and their corresponding dimensionally
reduced functions $u : Q \rightarrow \mathbb{R}$.
To this end,
we define the image of $H^1_{\textup{sym}}(\widehat{Q})$ under the dimensional
reduction as $V$.

Due to the coordinate transformation it is natural to endow  $V$ with the
pullback norm. This leads to the weighted inner product defined by
$(v,w)_r = \int_\omega vwr\, \text{d}(r,z)$ and the induced norm
$\Vert v \Vert_{L^2_r(Q)} = \sqrt{(v,v)_r}$. 
With this norm, we now define the space 
\begin{equation}
H^1_r(Q) = \{ u: Q\rightarrow \mathbb{R} \text{ is weakly differentiable and } \Vert u \Vert_{H^1_r} < \infty \}
\end{equation}
with the norm $\Vert v \Vert^2_{H^1_r(Q)} = \Vert v \Vert^2_{L^2_r(Q)} + \Vert \nabla v \Vert^2_{L^2_r(Q,\mathbb{R}^d)} $. 
Due to the weak differentiability of the reduced functions and the definition of the weighted norm, we have $V\subset H^1_r(Q)$. Our goal now is to show that we can identify this space with $H^1_\text{sym}(\widehat{Q})$, i.e. that $V = H^1_r(Q)$.  
For this we show, that for every function $u \in H^1_r(Q)$ its rotational extension $\widehat{u}$ defined by
\begin{equation}
\widehat{u}(x_1,x_2,x_3) = u(\vert (x_1,x_2) \vert, x_3) \text{ for } (x_1,x_2,x_3) \in \widehat{Q}
\end{equation}
belongs to $H^1(\widehat{Q})$. Due to the construction of the weighted norm it is only left to show that $\widehat{u}$ is also weakly differentiable. 

We define $\widehat{Q}_\varepsilon = \widehat{Q} \cap \{ x \in \widehat{Q}: \sqrt{x_1 +x_2} > \varepsilon \}.$ For $\varepsilon >0$ the coordinate transformation $\Phi$ restricted to $\Phi^{-1}(\widehat{Q}_\varepsilon)$ is differentiable, therefore $\widehat{u}$ is weakly differentiable on $\widehat{Q}_\varepsilon$. We have for a test function $\widehat{\varphi} \in C^\infty_c(\widehat{Q},\mathbb{R}^3)$
\begin{align*}
\int_{\widehat{Q}} \widehat{u}\text{ div}\widehat{\varphi}  \, \text{d}x & = \int_{\widehat{Q}\backslash \widehat{Q}_\varepsilon} \hat{u}\text{ div}\widehat{\varphi}  \, \text{d}x - \int_{\widehat{Q}_\varepsilon} \nabla \hat{u}_\varepsilon \cdot \widehat{\varphi}  \, \text{d}x + \int_{\Gamma_\varepsilon} \widehat{u} \widehat{\varphi}\cdot n \, \text{d}s,
\end{align*}
with $\Gamma_\varepsilon = \partial \widehat{Q}_\varepsilon \backslash \partial \widehat{Q}$ and $\nabla \widehat{u}_\varepsilon$ the weak gradient of the $\hat{u}$ restricted to $\widehat{Q}_\varepsilon$. Since $\widehat{u} \in L^2(\widehat{Q})$ and  $\widehat{\varphi} \in C^\infty_c(\widehat{Q},\mathbb{R}^3)$ the first term vanishes as $\varepsilon \rightarrow 0$.  We then define $\nabla \widehat{u}$ as the weak limit of $\nabla \widehat{u}_\varepsilon$. We claim that $\nabla \widehat{u}$ is the weak derivative  of $\widehat{u}$. This is the case if the boundary term vanishes as $\varepsilon \rightarrow 0$. 

To show this we use that $\Gamma_\varepsilon$ is a surface of revolution to deduce that for a function $\widehat{\psi}\in H^1_{\text{sym}}(\widehat{Q})$ we have $\int_{\Gamma_\varepsilon} \widehat{\psi}\ \cdot n \, \text{d}s  = 0$. We can then derive the following estimate:
\begin{equation}\label{eq:estimate_bdy}
\begin{aligned}
\left( \int_{\Gamma_\varepsilon} \widehat{u}\widehat{\varphi} \cdot n \, \text{d} s \right)  ^2 &= \left( \int_{\Gamma_\varepsilon} \widehat{u} (\widehat{\varphi} - \widehat{\varphi} (0,z))\cdot n \, \text{d} s \right)^2 \\& \le \int_{\Gamma_\varepsilon} \widehat{u}^2\, \text{d} s \int_{\Gamma_\varepsilon} ( (\widehat{\varphi} - \widehat{\varphi} (0,z))\cdot n )^2\, \text{d} s \\& \le \int_{\Gamma_\varepsilon} \widehat{u}^2 \, \text{d}s \int_{\Gamma_\varepsilon}(\varepsilon \Vert \nabla \widehat{\varphi}  \Vert_{L^\infty(\widehat{Q})} )^2 \, \text{d}s \\ & \le c(\widehat{\varphi} ) \varepsilon^2 \int_{\Gamma_\varepsilon}1 \, \text{d}s \int_{\Gamma_\varepsilon} \widehat{u}^2 \, \text{d}s  \le c(\widehat{\varphi} ,\widehat{Q}) \varepsilon^3 \int_{\Gamma_\varepsilon} \widehat{u}^2 \, \text{d}s
\end{aligned}
\end{equation}
since $\widehat{\varphi}  \in C^\infty(\widehat{Q},\mathbb{R}^3)$.

It can be checked that the constants appearing in the trace inequality for the boundary $\Gamma_\varepsilon$ depend on the parameter $\varepsilon^{-1}$, i.e.  
\begin{equation*}
\Vert \widehat{u} \Vert^2_{L^2(\Gamma_\varepsilon)} \le c\varepsilon^{-1} \Vert \widehat{u} \Vert^2_{H^1(\widehat{Q} \backslash \widehat{Q}_\varepsilon)}.
\end{equation*} 
This follows by deriving the trace estimates with regard to the weighted norms, which involves a derivative of the factor $r$, so that an upper bound for $r^{-1}$ needs to be estimated.

With this,  we deduce from the estimate \eqref{eq:estimate_bdy} that the boundary term $\left( \int_{\Gamma_\varepsilon} \widehat{u}\widehat{\varphi} \cdot n \, \text{d} s \right)  ^2$ vanishes as $\varepsilon \rightarrow 0$.

This means, that
for every function $u \in H^1_r(Q)$  the corresponding rotated function
$\widehat{u}: \widehat{Q} \rightarrow \mathbb{R}$ satisfies
$\widehat{u} \in H^1_{\textup{sym}}(\widehat{Q})$, such that
\begin{equation*}
\Vert \widehat{u} \Vert^2_{H^1(\widehat{Q})} = 2\pi\Vert u\Vert_{H^1_r(Q)}^2.
\end{equation*}

For rotationally symmetric sub-domains $\Omega \subset \widehat{Q}$ we
denote the transformed and dimensionally reduced domain with
$\omega \subset Q$,
which is one half of the cross section. The domain $\omega$ now has the boundary
$\partial \omega = \Gamma_{\textup{axis}} \cup \Gamma_{\textup{out}}$, where
$\Gamma_{\textup{axis}}$ corresponds to the axis of rotation and
$\Gamma_{\textup{out}}$ to the transformed boundary of the initial domain. In reverse, for a domain $\omega \subset Q$, we will denote its corresponding rotated three dimensional domain by $R(\omega) \subset \widehat{Q}$.
 
Lastly we shortly comment on the weak formulations of the reduced state equations. In particular for the Poisson problem
\begin{equation}\label{eq:pp_3d}
-\Delta \widehat{u} = f \text{ in } \Omega, \quad
\widehat{u} = 0 \text{ on } \partial \Omega
\end{equation}
the reduced formulation is given by
\begin{equation}\label{eq:pp_2d}
-\text{div}(r\nabla u)= r f \text{ in } \omega, \quad
u= 0 \text{ on } \Gamma_{\textup{out}}.
\end{equation}
This leads to the weak formulation for which a rotationally symmetric solution $u\in H^1_r(\omega)$ solves
\begin{align*}
\int_\omega \nabla u \cdot \nabla \phi r \, \text{d}(r,z) = -\int_{\omega} \text{div}(r\nabla u)\phi \, \text{d}(r,z)= \int_\omega f\phi r \, \text{d}(r,z)
\end{align*}
for all test functions $\phi \in C^1_{\Gamma_{\text{out}}}(\omega)$. In particular, no boundary condition arises on $\Gamma_{\text{axis}}$.

\section{Existence and Numerical Approximation of Optimal Domains}
\label{chap:ex_theory}

For an eigenvalue problem arising in a model of optimal insulation, we
now discuss how to establish existence of an optimal domain.

We look at the non-linear eigenvalue problem arising in optimal insulation and follow \cite{BBN17} closely for this section. We try to surround a
heat conducting body with an insulating material to get the best insulating
properties, i.e. to minimize the heat decay rate, which is given for the thickness of the insulating layer $\ell:
\partial \Omega \rightarrow \mathbb{R}_{+}$ by the
principal eigenvalue of the corresponding differential operator
\begin{equation}\label{eq:rayleigh_oi_l}
\lambda_\ell = \inf \left\lbrace  \int_\Omega \vert \nabla u \vert ^2 \  \, \text{d}x  + \int_{\partial \Omega} \ell^{-1}u^2 \, \text{d}s \ : \ \int_\Omega u^2 \, \text{d}x = 1 \right\rbrace.
\end{equation} 
The boundary term corresponds to Robin-type boundary conditions which result
from a model reduction in which the thickness $\ell$ with total mass $m$ is proportional
to the heat flux through the boundary.
With Hölder's inequality we can see that for a fixed $ u \in H^1(\Omega)$ 
the optimal thickness $\ell$ is given by 
\begin{equation*}
\ell(z) = \frac{m\vert u (z) \vert}{\int_{\partial \Omega} \vert u \vert \, \text{d}s}.
\end{equation*}
Thus, the optimal insulation can be obtained from a solution of the eigenvalue
problem
\begin{equation} \label{eq:rayleigh_oi}
\lambda_m = \min \left\lbrace J_m(u) = \int_{\Omega} \vert \nabla u\vert^2 \, \text{d}x + \frac{1}{m}\Vert u \Vert_{L^1(\partial \Omega)}^2 \ : \ \int_{\Omega} \vert u \vert ^2 \, \text{d}x = 1  \right\rbrace.
\end{equation}
We note, that the eigenfunction $u$ can be chosen to be non-negative. The
existence of this eigenfunction follows with the direct method of the calculus
of variations.
\begin{remark}\rm
With the transformation formula we can infer the following scaling property for
the eigenvalue. For $t > 0$
\begin{equation*}
t^{-2}\lambda_m(\Omega) = \lambda_{mt^d}(t\Omega).
\end{equation*}
This is the same scaling property as known from the eigenvalues of the
Dirichlet Laplacian or the Neumann Laplacian (c.f. \cite{henrot}), as long as
the mass of insulating material is scaled accordingly.
\end{remark}

Before proving existence and deriving the dimensionally reduced problem, we remark on the rotational symmetry of the eigenfunction $u$. In
\cite{BBN17} it was proven, that for a ball and for $m$ small enough, the
eigenfunction is not radial. However, experiments in \cite{BB19} indicate
that a rotationally symmetric solution exists. We adapt the optimization problem to only search for an eigenfunction among rotationally symmetric
functions, i.e. we look at the minimization problem 
\begin{align*}
\lambda^{\text{sym}}_m = \min \Big\lbrace  J_m(u) =  \int_{\Omega} \vert \nabla u\vert^2 \textup{ d}x + \frac{1}{m}\Vert u \Vert_{L^1(\partial \Omega)}^2: &\int_\Omega \vert u\vert^2\, \text{d}x =1,  \\  & u \text{ rotationally symmetric} \Big\rbrace.
\end{align*} 
This restriction may lead to larger eigenvalues and therefore to a larger optimal value for the shape optimization problem. However, even with the additional
constraint, the numerical results of the dimensionally reduced problem have been
consistent with the results we expect from the three-dimensional shape
optimization problem, see Section \ref{chap:c_num_ex}.

The corresponding shape optimization problem for a fixed mass $m > 0$ is
defined as follows:
\begin{align*} \label{problem_m}
\textup{Minimize } &   \lambda^{\text{sym}}_m( \Omega) = J_m(u, \Omega) = \int_{ \Omega} \vert \nabla u \vert^2 \, \text{d}x + \frac{1}{m} \left( \int_{\partial  \Omega} \vert u \vert \, \text{d}s \right)^2 \\  \tag{$\mathbf{\widehat{P}_m}$}
\textup{w.r.t } &  \Omega \subset \widehat{Q} \subset \mathbb{R}^3 \text{ open, convex and rotationally symmetric} \\ \text{and } & \vert \Omega \vert = M \\   \text{ s.t.  }& u \in H_{\text{sym}}^1( \Omega) \text{ is an eigenfunction to } \lambda^{\text{sym}}_m(\Omega) \text{ with } \Vert u \Vert_{L^2(\Omega)} = 1
\end{align*}

Here $\widehat{Q} \subset \mathbb{R}^3$ is an open, rotationally symmetric and
bounded hold-all domain.
The condition that $u$ is an eigenfunction  is equivalent to the minimality of  $u \in H_{\text{sym}}^1(\Omega)$ with $\Vert u \Vert_{L^2(\Omega)} =1$ for $J_m(\ \cdot, \Omega)$ for a fixed domain $\Omega.$ 

To prove existence, we adapt the strategy from \cite{BW18}. However,
due to the lack of homogeneous Dirichlet boundary conditions, which allow for trivial extensions in $H^1(\widehat{Q})$, we need to
incorporate a convergence result for special functions of bounded variations.
This approach is often used for eigenvalue problems with a Robin-type
boundary condition (see e.g. \cite{BG}), since the boundary term occurring in
the eigenvalue problem allows for the use of the compactness results of $SBV$. 

\begin{proposition} \label{thm:oi_ex}
There exists an optimal pair $(\Omega,u)$ for \eqref{problem_m}.
\end{proposition}
\begin{proof}
We can select a minimizing sequence $(\Omega_n,u_n)_{n \in \mathbb{N}}$ of convex domains $\Omega_n$ and eigenfunctions $u_n \in H_{\text{sym}}^1(\Omega_n)$ with $\Vert u_n \Vert_{L^2(\Omega_n)}=1$ for $n\in \mathbb{N}$. After passing to a subsequence we find an open, convex and rotationally
symmetric domain $\Omega \subset \widehat{Q}$, such that
$\chi_{\Omega_n} \rightarrow \chi_{\Omega}$ in $L^1(\widehat{Q})$,
see \cite[Lemma 3.1]{BG97}. Therefore we also
maintain the volume $\vert \Omega \vert =M$.
 Furthermore, we use that  we can chose $u_n$ to be non-negative.
After trivially extending $u_n \in H_{\text{sym}}^1(\Omega_n)$ to $\tilde{u}_n \in SBV(\widehat{Q})$,
we have for all $n$, we can find a suitable bound $C  < \infty$ such that
\begin{align*}
\Vert \widetilde{u}_n \Vert_{L^2(\widehat{Q}) }=1,\\ 
\Vert \nabla \widetilde{u}_n \Vert_{L^2(\widehat{Q}, \mathbb{R}^3)} &\le C.
\end{align*}
  Here, $\nabla \widetilde{u}$ refers to the piecewise weak gradient rather than the
 weak gradient.
 From \cite[Theorem 2.6]{BG97} we can deduce that the measure $\vert D \chi_{\Omega_n} \vert$ coincides with $\mathcal{H}^{N-1} \llcorner \partial \Omega_n$. Since the functions $\widetilde{u}_n$ are weakly differentiable on $\Omega_n$ and $\widehat{Q} \backslash \overline{\Omega_n}$, we can therefore identify the jump set of $\widetilde{u}_n$ with the boundary of $\Omega_n$. The eigenfunctions $u_n$ are chosen to be non-negative and $\widetilde{u}_n = 0$ on $\widehat{Q}\backslash \overline{\Omega_n}$. Since $u_n$ is a minimizing sequence of eigenfunction, we can then bound the boundary terms
\begin{equation*}
\int_{J_{\widetilde{u}_n}}\ \widetilde{u}_n^+\nu^+ - \widetilde{u}_n^-\nu^-\, \text{d}s = \int_{\partial \Omega_n} \widetilde{u}_n \, \text{d}s =  \int_{\partial \Omega_n} \vert \widetilde{u}_n \vert \, \text{d}s \le \sqrt{m J_m(u_n,\Omega_n)} \le C
\end{equation*}
 for the unit normals $\nu^+,\nu^- $ along the jump sets $J_{\widetilde{u}_n}$, see e.g. \cite[Example 10.2.1]{attouch}.

Since $Dv(\widehat{Q}) = \int_{\widehat{Q}} \nabla v\, \text{d}x + \int_{J_v\cap \widehat{Q}} v^+\nu^+-v^- \nu^-\, \text{d}s $
for all $v \in SBV(\widehat{Q})$, the sequence $(\widetilde{u}_n)_{n\in \mathbb{N}}$
is bounded in $SBV(\widehat{Q})$, so that we can use the compactness theorem
for special functions of bounded variations \cite[Theorem 2.1]{BG} to find
a function $\widetilde{u} \in SBV(\widehat{Q})$, s.t.
\begin{align}
D\widetilde{u}_n \rightharpoonup^\star D\widetilde{u} \quad &\text{ in the sense of measures}  \label{eq:oi_conv_meas}\\
\widetilde{u}_{n_k} \rightarrow \widetilde{u} \quad &\text{ strongly in } L^1(\widehat{Q})   \label{eq:oi_l1_conv}\\
\nabla \widetilde{u}_{n_k} \rightharpoonup \nabla \widetilde{u} \quad &\text{ weakly in } L^2(\widehat{Q}, \mathbb{R}^N) \label{eq:oi_conv_grad}\\
\mathcal{H}^{N-1} (J_{\widetilde{u}} ) \le  & \liminf_{n \rightarrow \infty} \mathcal{H}^{N-1} (J_{\widetilde{u}_{n_k}} ).
\end{align}
Due to the boundedness of $\widetilde{u}_n$ with respect to the $L^2$-norm, we
further have that
\begin{equation}\label{conv_l2_oi}
\widetilde{u}_{n_k} \rightharpoonup \widetilde{u} \text{ weakly in } L^2(\widehat{Q}).
\end{equation}
This implies with $\chi_{\Omega_n} \rightarrow \chi_{\Omega}$ in $L^1(\widehat{Q})$
and \eqref{eq:oi_l1_conv}, that $\widetilde{u}\vert_{\widehat{Q} \backslash
\overline{\Omega}} = 0$.

Next we show, that $u:= \widetilde{u}\vert_{\Omega} \in H^1(\Omega)$.
Let $\phi$ be a test function from $C^\infty_c(\Omega)$. Then for  all
$n \ge N$ with $N$ sufficiently large we have, due to the convexity,
that $\phi \in C^\infty_c(\Omega_n)$ (c.f. \cite[Lemma 4.2]{BG97}) and
\begin{align*}
\int_{\Omega}  u\ \text{div} \phi  \, \text{d}x &= \int_{\widehat{Q}}  \widetilde{u} \ \text{div}\phi \, \text{d}x = \lim_{n \rightarrow \infty} \int_{\widehat{Q}} \widetilde{u}_n\ \text{div}\phi \, \text{d}x = \lim_{n \rightarrow \infty} \int_{\Omega_n} \widetilde{u}_n \ \text{div} \phi \, \text{d}x \\ &= \lim_{n \rightarrow \infty} -\int_{\Omega_n} \nabla \widetilde{u}_n \cdot  \phi \, \text{d}x = \lim_{n \rightarrow \infty} - \int_{\widehat{Q}} \nabla \widetilde{u}_n  \cdot  \phi \, \text{d}x =  - \int_{\Omega} \nabla u \cdot \phi \, \text{d}x,
\end{align*}
i.e. the weak gradient coincides with the approximate gradient on $\Omega$.

The rotational symmetry of the eigenfunctions $(u_n)_{n\in\mathbb{N}}$ is preserved under $L^1$-convergence, and therefore $u \in H_{\text{sym}}^1(\Omega )$.

Because $u \in H^1(\Omega)$ and $\widetilde{u}\vert_{\widehat{Q} \backslash \overline{\Omega}} =0$,
we have that $J_{\widetilde{u}} \subset \partial \Omega$. Because of $\widetilde{u}\vert_{\widehat{Q} \backslash \overline{\Omega}} =0$,
we have for the trace on $\partial \Omega \backslash J_{\widetilde{u}}$ that $\widetilde{u} = 0$.
This results in $\int_{J_{\widetilde{u}}} \widetilde{u} \, \text{d} \mathcal{H}^{d-1} = \int_{\partial \Omega} u \, \text{d} \mathcal{H}^{d-1}$.
Then \eqref{eq:oi_conv_meas} and \eqref{eq:oi_conv_grad}, and \cite[Theorem 2.6]{BG97} imply
\begin{equation*}
\int_{\partial \Omega_n} u_n \, \text{d}s = \int_{J_{u_n}} u_n \, \text{d} \mathcal{H}^{d-1} \rightarrow \int_{J_{u}} u \, \text{d} \mathcal{H}^{d-1}.
\end{equation*}
By the assumption that the eigenfunctions $u_n$ are non-negative this means that
\begin{equation} \label{conv_boundary_oi}
\Vert u_n \Vert_{L^1(\partial \Omega_n)} \rightarrow \Vert u \Vert_{L^1(\partial \Omega)}.
 \end{equation}

To show that $u \neq 0$, we follow an argument in \cite[Proposition 1]{BG10}
and show that $u_n \rightarrow u$ in $L^2(\widehat{Q})$.
We note that for the minimizing sequence $(u_n)_{n\in \mathbb{N}}$  we have
that $u_n^2 \in SBV(\widehat{Q})$. This is due to $\Vert u_n \Vert_{L^2(\widehat{Q})} = 1$ and the total variation
\begin{align*}
 D(u_n^2)(\widehat{Q}) &= \int_{\widehat{Q}} 2 u_n \nabla u_n \, \text{d}x + \int_{J_{u_n}\cap \widehat{Q}}u^2_n\, \text{d}s \\ & \le c\left(\int_{\widehat{Q}} \vert \nabla u_n \vert ^2 \, \text{d}x +\int_{\widehat{Q}} u^2 \, \text{d}x \right) + \int_{J_{u_n} \cap \widehat{Q}}u^2_n\, \text{d}s.
\end{align*}

Using results from \cite[Equations (1.5) and (1.6)]{payne}, we can bound the constant of the trace
inequality for the functions $(u_n)_{n \in \mathbb{N}}$  independently of
$\Omega_n$, so that 
\begin{equation}\label{payne}
\int_{\partial \Omega_n} u^2_n\, \text{d}s \le C(\widehat{Q}) \Vert u_n \Vert^2_{H^1(\widehat{Q})}.
\end{equation}
Due to geometric constraints of the domains the constant $C(\widehat{Q})$ can indeed be chosen independently of $\Omega_n$: In \cite{payne} the divergence theorem is used for a fixed convex domain $\Omega^\prime$ with the vector field $f_{\Omega^\prime}(x)=x-x_{\Omega^\prime}$  for a point $x_{\Omega^\prime}\in \Omega^\prime$. For this vector field we have $f_{\Omega^\prime}.\nu \ge k(\Omega^\prime) >0$ a.e. on $\partial \Omega^\prime$, where $\nu$  is the outer unit normal vector on $\partial \Omega^\prime$. 
 Using the convexity, boundedness and fixed volume of the admissible domains, we can find uniform bounds on the radius of an incircle and the diameter, c.f. the Steinhagen inequality, \cite{steinhagen}, and \cite[Theorem 50]{eggleston_1958}.
 Thus, we can choose $x_{\Omega'}$ as a center of an incircle. Consequently, the mentioned bounds can be used to get a lower bound for $k(\Omega')$ which is independent of $\Omega'$ and this leads to the constant $C(\hat Q)$. \\
With this trace estimates, since since $(u_n)_{n \in \mathbb{N}}$ is a minimizing  sequence, the sequence
$(u^2_n)_{n \in \mathbb{N}}$ is bounded in $BV(\widehat{Q})$ and admits a subsequence,
which converges weakly to a function $v$ in $BV(\widehat{Q})$. Especially, since $\widehat{Q}$ is bounded, due to the compact embedding of $BV(\widehat{Q})$ into $L^1(\widehat{Q})$, we have $u_n^2 \rightarrow v$ in $L^1(\widehat{Q})$. Due to the assumed non-negativity of $u$, this results in $u_n \rightarrow \sqrt{v}$ strongly in $L^2(\widehat{Q})$. With \eqref{conv_l2_oi} and the uniqueness of the limit, this results in the strong convergence of $u_n \rightarrow u$ in $L^2(\widehat{Q})$. 

We can now use \eqref{eq:oi_conv_grad}, \eqref{conv_l2_oi} and
\eqref{conv_boundary_oi} to show, that $u$ satisfies the variational eigenvalue
equation for the eigenvalue 
\begin{equation}\label{eq:ev_oi_liminf}
\lambda_m^{\text{sym}}(\Omega) = \liminf_{n \rightarrow \infty} \lambda_m^{\text{sym}}(\Omega_n),
\end{equation} which proves, that $(\Omega,u)$ is an admissible pair. The optimality of the pair follows from \eqref{eq:ev_oi_liminf}, since $(\Omega_n,u_n)_{n \in \mathbb{N}}$ was chosen as an infimizing sequence.
\end{proof}

\section{Discretized  Reduced Problem}
\label{chap:dis_red_p}

Next, we derive the dimensionally reduced problem and define the numerical scheme and point out technical difficulties in stability. Lastly, we address how this scheme can be applied to other optimization problems

\subsection{Dimensionally Reduced Problem}
\label{sec:dim_red}

After transformation and dimensional reduction, we obtain the equivalent minimization problem
\begin{align*}\label{problem_m_red}
	\textup{Minimize } &   \lambda^r_m(\omega) = J_m^r(u,\omega) = \int_{\omega} \vert \nabla u \vert^2 r \, \text{d}(r,z) + \frac{2\pi}{m} \left( \int_{\Gamma_{\textup{out}}} \vert u \vert r \, \text{d}s \right)^2 \\ 
\textup{w.r.t } & \omega \subset Q\subset \mathbb{R}_+\times\mathbb{R} \text{ open and convex} \\ \text{and } & 2\pi\vert \omega \vert_r = 2\pi \int_\omega r \, \text{d}(r,z)= M \tag{$\mathbf{P_m}$}   \\ 
 \text{s.t. } &u \in H^1_r(\omega) \text{ is an eigenfunction to } \lambda^r_m \text{ with } \int_\omega \vert u \vert^2 r \, \text{d}(r,z) = 1 \\ \text{ and } & \text{the rotated domain  } R(\omega) \subset  \mathbb{R}^3 \text{ is also convex}.
\end{align*}

We note, that for the distribution $\ell$ of insulating material we now have
\begin{equation*}
	2\pi\int_{\Gamma_{\textup{out}}} \ell r \, \text{d}s = m \quad  \text{and} \quad
	\ell(z) = \frac{m \vert u \vert}{2\pi\int_{\Gamma_{\textup{out}}} \vert u \vert r \, \text{d}s }.
\end{equation*}
We introduce a regularization for numerical treatment, c.f. \cite{BB19}, and for $\varepsilon >0$  we look for a minimizer  $u \in H^1_r(\omega)$ with $\Vert u \Vert_{L^2_r(\omega)} = 1$  of the differentiable functional
\begin{equation*}
	J^r_{m,\varepsilon}(u) = \Vert \nabla u \Vert_r^2 + \frac{2\pi}{m}   \Vert u \Vert^2_{r,L^1_\varepsilon(\Gamma_{\textup{out}})}
\end{equation*}
with the regularized norm
\begin{equation*}
	\Vert u \Vert_{r,L^1_\varepsilon(\Gamma_{\textup{out}})} = \int_{\Gamma_{\textup{out}}} \vert u \vert_\varepsilon r \, \text{d}s \quad \text{ with }\vert u \vert_\varepsilon = \sqrt{u^2 + \varepsilon^2}.
\end{equation*}

A minimizer satisfies for all $v \in H^1_r(\omega)$ the variational formulation
\begin{equation}\label{eq:reg_var}
	\left( \nabla u, \nabla v\right)_r + \frac{2\pi}{m}\Vert u \Vert_{r,L^1_\varepsilon(\Gamma_{\textup{out}})} \int_{\Gamma_{\textup{out}}} \frac{uv}{\vert u \vert_\varepsilon} r \, \text{d}s = \lambda^{r,\varepsilon}_m(\omega) \left(u,v\right)_r.
\end{equation}
As $\varepsilon \rightarrow 0$ the eigenvalue $\lambda_m^r$ is approximated using a gradientflow to find a function $u \in H^1_r(\omega)$ solving \eqref{eq:reg_var}.

This leads to a regularized version of the shape optimization problem \eqref{problem_m_red} depending on $\varepsilon >0$, which will be discretized in the next section.
The effect of the regularization can be controlled by using the unconditional uniform estimate $0 \le \vert u \vert_{\varepsilon} - \vert u \vert \le \varepsilon$. For more details on the iterative minimization and discretization we refer to \cite{BB19}, since the results carry over to the reduced problem. We will not go into further detail
here and only mention what is necessary to define the discrete shape optimization scheme and discuss aspects of stability of the discretization and the iterative approximation of the optimal domain.

\subsection{Spatial Discretization}

Following \cite{BB19} we approximate $\omega$ with a polyhedral domain $\omega_h$ and, given a regular
triangulation $\mathcal{T}_h$, we define the finite element space
\begin{equation*}
\mathcal{S}^1(\mathcal{T}_h) = \left\lbrace v_h \in C(\overline{\omega}_h) : v_h \vert _T \in P_1(T) \text{ for all } T \in \mathcal{T}_h \right\rbrace.
\end{equation*}
Including a quadrature formula we consider the functional
\begin{equation*}
	J^r_{m,\varepsilon,h}(u_h) = \Vert \nabla u_h \Vert_{L^2_r(\omega_h)}^2 + \frac{2\pi}{m} \Vert u_h \Vert^2_{r,L^1_{\varepsilon,h}(\Gamma_{\textup{out},h})}
\end{equation*}
with the discretized and regularized $L^1$-norm
\begin{equation*}
	\Vert u_h \Vert_{r,L^1_{\varepsilon,h}(\Gamma_{\textup{out},h})} = \int_{\Gamma_{\textup{out},h} } \mathcal{I}_h \vert u_h \vert_\varepsilon r \, \text{d}s = \sum_{z \in \mathcal{N}_h \cap \Gamma_{\textup{out},h}} \beta_z \vert u_h(z) \vert_\varepsilon
\end{equation*}
with the nodal interpolation operator $\mathcal{I}_h: C(\overline{\omega}_h)\rightarrow \mathcal{S}^1(\mathcal{T}_h)$
corresponding to the nodal basis functions $\varphi_z \in \mathcal{S}^1(\mathcal{T}_h)$ and $\beta_z := \int_{\Gamma_{\textup{out},h}} \varphi_z r \, \text{d} s$.
The corresponding variational formulation is given by
\begin{equation}\label{eq:reg_var_dis}
	\left( \nabla u_h, \nabla v_h\right)_r + \frac{2\pi}{m}\Vert u_h \Vert_{r,L^1_{\varepsilon,h}(\Gamma_{\textup{out},h})} \int_{\Gamma_{\textup{out},h}} \frac{u_hv_h}{\vert u_h \vert_\varepsilon} r \, \text{d}s = \lambda^{r,\varepsilon,h}_m(\omega_h) \left(u_h,v_h\right)_r
\end{equation}
for all $v_h,u_h \in \mathcal{S}^1(\mathcal{T}_h)$.
\color{black}

We can now define the discretized shape optimization problem:
\begin{align*}\label{problem_m_discrete}
\textup{Minimize } &  \lambda^{r,\varepsilon,h}_m(\omega_h) \\
\textup{w.r.t. }&\omega_h \subset \mathbb{R}_+\times \mathbb{R}, \ \mathcal{T}_h \in \mathbb{T}_{c_{\textup{usr}}} \textup{ triangulation of } \omega_h \tag{$\mathbf{P_{m,h,\varepsilon}}$} \\
\textup{s.t. }& u_h  \in \mathcal{S}^1(\mathcal{T}_h) \textup{ solves } \eqref{eq:reg_var_dis} \textup{ and } \Vert u_h \Vert_{L^2_r(Q)} = 1
 \\ & \omega_h \subset Q \textup{ is  convex and open and } 2\pi\vert \omega_h \vert_r = M \\ \text{and } & \text{the rotated domain  } R(\omega_h) \subset  \mathbb{R}^3 \text{ is also convex.} 
\end{align*}
Here, $\mathbb{T}_{c_{\textup{usr}}}$ is the class of conforming, uniformly
shape regular triangulations $\mathcal{T}_h$ of polyhedral subsets of
$\mathbb{R}^2$ with $h_T/\varrho_T \le c_{\textup{usr}}$  for all
elements $ T \in \mathcal{T}_h$ with diameter $h_T \le h$ and inner radius $\varrho_T$ for a universal
constant $c_{\textup{usr}} >0$.

We adapt the numerical approximation of the eigenvalue arising in optimal insulation from \cite{BB19} to the dimensional reduced eigenvalue problem. However, the dimensional reduction makes it difficult to infer the consistency and stability results for the dimensionally reduced eigenvalue problem and the shape optimization problem. \\
While we are able to estimate the interpolation error in the weighted norm with the
interpolation error regarding the $H^1$-norm, see e.g. \cite[Theorem 3.2]{pde} for functions in $H^1(\omega)$, this does not provide a sufficient result for functions in $H^1_r(\omega)$. \\
The lack of an error estimate using the weighted norm, which is needed for the  $\Gamma$-convergence of the discrete functionals, c.f. \cite[Corollary 4.2]{BB19},  poses additional difficulties for the convergence analysis here. 

There are some results for interpolation estimates regarding weighted norms, such as \cite{elhatri},\cite{atamni},\cite{nochetto} or \cite{antil18}. From \cite[Theorem 4.1]{elhatri} for example we can derive  for a  uniformly shape regular family of triangulations $\mathcal{T}_h $ the estimate
\begin{equation*}
\vert v - \mathcal{I}_h v\vert_{H^m_r(K)} \le Ch^{2-m}\vert v \vert_{H^{2}_r(K)}
\end{equation*}
for all $v \in H^{2}_r(K)$, and a triangle $K$ in $\mathcal{T}_h $ of $\omega_h$ and the nodal interpolation operator $ \mathcal{I}_h$. However, this estimate is not sufficient to get the corresponding results with respect to the weighted norm, since it provides no estimates for the interpolation for the trace with respect to the weighted norm.

\subsection{Application to other shape optimization problems}

The shape optimization problem as described in the previous sections can be
applied to other suitably posed problems of the form \eqref{problem}. For a
minimizing sequence $(\Omega_n, u_n)$, the sequence of trivially extended
functions $ {u}_n$ should be bounded in $SBV(\widehat{Q}$) or $H^1_0(\Omega_n)$. In
order to use the compactness results of $SBV$, the jumps of the minimizing
functions need to be controlled. In the eigenvalue problem for optimal
insulation this condition is satisfied due to the boundary term occurring in
the eigenvalue which we want to minimize. However, the results from
\cite{payne} as used to obtain the bound on the trace \eqref{payne}, also
guarantee that the $BV$-norm is bounded. This means, rather than just using it
to prove the strong $L^2$-convergence, it also allows us to obtain a convergent
subsequence for shape optimization problems in which we have neither a boundary
term in the objective value nor a homogeneous Dirichlet boundary condition.

Furthermore, to guarantee existence of an optimal domain, we need suitable continuity of the state operator, such that an accumulation pair $(\Omega,u)$ of a minimizing sequence, $u \in H^1(\Omega)$
also solves the state equation in $\Omega$.

The objective functional has to be (weakly) lower semi-continuous (depending on
the mode of convergence of $u_n \rightarrow u$) to guarantee optimality of the
limit. Lastly, the consistency and numerical stability of the discrete scheme
has to be guaranteed, for example via strong continuity properties and density
results.

We will see in the next section that the shape optimization algorithm works
independently of the optimization problem itself, i.e. only the objective value
and the state equation need to be implemented specific to the optimization
problem.

\section{Iterative Computation of Optimal Domains}
\label{chap:num_approx}

We next address the iterative numerical approximation of optimal domains. After
dimensional reduction and spatial discretization, we obtain the following class of shape
optimization problems.
\begin{align*}
	\label{eq:P_h}
\textup{Minimize } & \int_{\omega_h} j_h((r,z),u_h(r,z),\nabla u_h(r,z)) r\, \text{d}(r,z) \\
\textup{w.r.t. }&\omega_h \subset \mathbb{R}_+\times\mathbb{R}, \ \mathcal{T}_h \in \mathbb{T}_{c_{\textup{usr}}} \textup{ triangulation of } \omega_h  \tag{$\mathbf{P_h}$}\\
\textup{s.t. }& u_h\in \mathcal{S}^1(\mathcal{T}_h) \textup{ solves the respective discrete state equation} \\  & \omega_h \subset Q \textup{ is  convex and open and } 2\pi\vert \omega_h \vert_r = M \\ \text{and } & \text{the rotated domain  } R(\omega_h) \subset  \mathbb{R}^3 \text{ is also convex,} 
\end{align*}
where $j_h$ now denotes the discrete transformed density function of \eqref{problem} and with
$\mathbb{T}_{c_{\textup{usr}}}$ the class of conforming, uniformly shape
regular triangulations $\mathcal{T}_h$ of polyhedral subsets of $\mathbb{R}^2$
with $h_T / \varrho_T \le c_{\textup{usr}}$ for all $T \in \mathcal{T}_h$ for a universal constant
$c_{\textup{usr}} >0$.

We adopt an approach  similar to \cite{BW18}, where the admissible domains
are obtained from a discrete deformation of a given convex reference domain. A
convex polygonal domain $\omega_h$ with a regular triangulation $\mathcal{T}_h $
is optimized by moving the vertices of the triangulation. For a piecewise
linear deformation field $V_h \in \mathcal{S}^1(\mathcal{T}_h)^2$ the
triangulation of the updated domain is obtained by a piecewise linear
perturbation $T_t = I + tV_h$ of the domain. The vertices of the updated
triangulation are given by $x_i + tV_h(x_i), i =1,\dots,N$.

Rather than deforming the entire triangulation, we deform the boundary
$\Gamma_{\textup{out}}$ of $\omega_h$, and then generate a triangulation of
$\omega_h$, to calculate the objective values or to find the deformation field.
Equivalently, we could also say that we add remeshing of the domain to the
deformed triangulations of \cite{BW18}. So rather than trying to solve
\eqref{eq:P_h}, we instead solve the problem as follows.
\begin{align*}
\text{Minimize } & \int_{\omega_h} j_h((r,z),u_h(r,z),\nabla u_h(r,z)) r\, \text{d}(r,z) \\ \text{w.r.t. } & \Phi_h \in  \mathcal{S}^1(\mathcal{T}_h(\widehat{\omega}))^2, \ \mathcal{T}_h(\omega_h) \in \mathbb{T}_{c_{\textup{usr}}} \text{ a triangulation of }\omega_h  \\
\text{s.t. } &\Vert D\Phi_h\Vert_{L^{\infty}(\widehat{\omega})} + \Vert [D\Phi_h]^{-1}\Vert_{L^{\infty}(\widehat{\omega})} \le c \\ & \omega_h = \Phi_h(\widehat{\omega}) \subset Q \textup{ is  convex and open and } 2\pi\vert \omega_h \vert_r = M \\ & \text{the rotated domain  } R(\omega_h) \subset  \mathbb{R}^3 \text{ is also convex}  \\
\text{ and } & u_h \in \mathcal{S}^1(\mathcal{T}_h(\omega_h)) \text{ solves the respective discrete state equation}.
\end{align*}
Here, $\mathcal{T}_h(\omega_h)$ and $\mathcal{T}_h(\widehat{\omega})$ are regular
triangulations generated to approximate $\omega$ and $\widehat{\omega}$. The triangulation $\mathcal{T}_h(\widehat{\omega})$ remains fixed during the optimization.
This  comes with a higher computational cost, due to the regular generation of the triangulation.
Since the deformation of the entire triangulation (as in \cite{BW18}) has
often led to a degeneration of the triangulation and the boundary nodes in the conducted experiments, a
frequent generation of a new triangulation was often necessary in either
versions.

The triangulation $\mathcal{T}_h(\omega_h)$
was generated by deforming a triangulation of the half-disk, since it allows for a good approximation of the boundary.

Since the approximation of the optimal domains with the described triangulation
seemed sufficient for the problems for which the optimal domain was already
known, only this approach was used. Whether this causes a geometric bias for
the approximated optimal domains was also not further investigated. The
solvability of the discretized shape optimization with deformed triangulations
is discussed in \cite{BW18}.

Furthermore, the boundedness of the admissible domains was also not included as a constraint in
the implemented code and we did not observe degeneration in the examples under consideration.

In the following sections we  look at the details of the optimization
algorithm. In Section \ref{sec:shape_gradient} we shortly introduce the notion
of shape gradients. After considering the convexity constraint in Section
\ref{sec:conv_constraint}, we look at how to find a suitable deformation field
in Section \ref{sec:def_field}. In Section \ref{sec:line_search} we state the necessary conditions which  determine the step size $\tau$ with which
to update the domain.

The implemented code was adapted from the algorithm described in
\cite{BW18} and the code used in \cite{BB19} for numerical
experiments. An implementable pseudo code is listed in Section \ref{sec:pseudo_code}.

\subsection{Shape Gradients} \label{sec:shape_gradient}

In order to find a suitable deformation field which leads to an optimized
domain, we first give a short summary of shape derivatives for the
PDE constrained shape optimization problems.

 The objective value of the minimization problem is given by the shape functional 
\begin{equation*}
J(\Omega) = \int_{\Omega} j(x,u(x),\nabla u(x)) \, \text{d}x
\end{equation*}
for a suitable cost function $j$ and the solution of the state equation $u \in H^1(\Omega)$.

Let $\Omega \subset \widehat{Q}$ be a fixed open and convex domain. Perturbations of identity $T_t = I + tV$ with $V \in C^{0,1}_c(\widehat{Q})$ lead to the
Eulerian derivative of the shape functional
\begin{equation*}
J'(\Omega,V) = \lim_{t \rightarrow 0} \frac{J(\Omega_t)- J(\Omega)}{t}
\end{equation*}
with the deformed domain $\Omega_t = T_t(\Omega)$ and 
\begin{equation*}
J(\Omega_t) = \int_{\Omega_t} j(x,u_t(x),\nabla u_t(x)) \, \text{d}x
\end{equation*}
where $u_t \in H^1(\Omega_t)$ is the solution of the state equation in
$\Omega_t$. The shape derivative can also be formulated in Hadamard form, i.e.
as a function on the boundary of the domain
\begin{equation*}
J'(\Omega,V) = \int_{\partial \Omega} g V.n \, \text{d}s
\end{equation*}
for an appropriate function $g$. 
The Hadamard derivative relies on certain regularity properties, but for finding a suitable
descent direction for our optimization problem this is neglected in our case. For more details on shape derivatives and shape sensitivity analysis we refer  to \cite{veodf} and \cite{zolesio}.

Both for the representation of the shape derivative on the volume and on the
boundary, the shape derivative is problem-specific. 
Therefore, we opt to only approximate the shape gradient on the boundary points
of $\Gamma_{\textup{out}}$ with a difference quotient. This involves a
high computational cost, but allows for 
different optimization problems to be approximated without having to adapt the
shape derivative. In some numerical experiments for the shape optimzation algorithm this approach has also led to better results in
optimization even for most of those problems, in which the Hadamard derivative
was beforehand known and could be approximated directly on the boundary.
The experiments documented in Section \ref{chap:c_num_ex} were also implemented so
that the shape gradient was approximated using forward algorithmic
differentiation, however without any notable difference in the approximated
optimal domains.

\subsection{Convexity Constraint}\label{sec:conv_constraint}

To ensure that the deformed domain $\omega_h$ is also convex, we need to
incorporate a constraint for the deformation field $V_h$. This approach follows
again \cite{BW18}.
Let $\omega_h \subset \mathbb{R}_+ \times\mathbb{R}$ be a simply connected
polygon and let $N$ be the number of boundary vertices of $\omega_h$ on
$\Gamma_{\textup{out}}$ with coordinates $x^i \in \mathbb{R}^2 , i = 1, \dots,
N$, in counter-clockwise order. It can be seen, that $\omega_h$ is convex if
and only if the interior angles are less than or equal to $\pi$.
 By using the cross product, this in turn is equivalent to
\begin{equation} \label{eq:char_conv_num}
C_i(X) := (x_1^{i-1} - x_1^i)(x_2^{i+1} - x_2^i)-(x_2^{i-1} - x_2^i)(x_1^{i+1} - x_1^i) \le 0
\end{equation}
for $i = 2, \dots, N-1$. 
 For the reduced optimization problems to be equivalent to the three-dimensional problems, we further need to guarantee that the corresponding three-dimensional rotated domain $R(\omega_h)$ is also convex. Therefore the interior angles for the nodes on the axis of rotation (i.e. where
$\Gamma_{\textup{out}}$ and $\Gamma_{\textup{axis}}$ intersect) have to be less than
or equal to $\pi/2$. This leads to the inequalities 
\begin{equation*}
C_1(X) := -2 x_1^2(x_2^2 - x_2^1)\le 0
\end{equation*}
for $i = 1$ and for $i = N$
\begin{equation*}
C_N(X) := 2 x_1^{N-1}(x_2^{N-1} - x_2^N)\le 0,
\end{equation*}
with the argument $X$ representing the vector $(x^1,\dots,x^N)$.

The last two inequalities are derived from \eqref{eq:char_conv_num} by using
the assumed symmetry of the corresponding three-dimensional domain.

 The convexity of the deformed domain $(I+t^0V_h) (\omega_h)$ is equivalent to
 $C_i(X + t^0V_h(X)) \le 0$. With a first-order expansion of this quadratic
 constraint we obtain the constraint
\begin{equation*}
C_i(X) + t^0DC_i(X)V_h(X) \le 0, \ \forall i = 1, \dots,N.
\end{equation*}
The constraint on the convexity of the three-dimensional domain  will be realized by having gliding boundary conditions on the nodes lying on the axis of rotation, so that they are only allowed to move along the axis of rotation and not away from it.
\\
However, for simplicity, the constraint that $R(\omega_h)$ is convex will be used in the definitions of the optimization problems, even if the convexity of the three-dimensional domain itself is not evaluated, only the conditions on the two-dimensional domain.

\subsection{Finding the Deformation Field}
\label{sec:def_field}

We follow \cite{BB19} to compute the deformation field $v$ from the linear functional $J'(\omega,\cdot)$. 
In order to satisfy a constraint on the volume of the three-dimensional domain, we
incorporate the constraint on the vector field, which relates to the
transformed divergence operator. So, rather than requiring $\text{div}(v) = 0$,
we instead search for deformation fields with $r^{-1}\text{div}(r v) = 0$. In order to satisfy the convexity constraint we have gliding boundary condition, i.e. we search for deformation fields $v \in H^1_{r,\text{glide}}(\omega)^d:=\{v \in H^1_r(\omega)^d: v_1 = 0 \text{ on } \Gamma_{\text{out}}\}$. 
 This means we find $v \in H^1_{r,\text{glide}}(\omega)^d $ and $q \in L^2_r(\omega)$ such that
\begin{align*}
\int_\omega v \cdot w \, \text{d}(r,z) + \int_\omega \nabla v : \nabla w \, \text{d}(r,z) - \int_\omega p \text{ div}(rw) \, \text{d}(r,z) &= -J'(\omega, w) \\ \int_\omega q \text{ div}(rv) \, \text{d}(r,z) &= 0
\end{align*}
for all $(w,q) \in H^1_{r,\text{glide}}(\omega)^d\times L^2_r(\omega)$.

Here, only the bilinear form which pertains to the divergence was transformed, since it turned out that
using the untransformed bilinear form provided better numerical results.

We discretize the system with the Crouzeix–Raviart method, i.e. we discretize
$L^2_r(\omega)$ with $\mathcal{L}^0(\mathcal{T}_h)$, the elementwise constant functions,
and $H^1_{r,D}(\omega)$ with the non-conforming space
\begin{align*}
\mathcal{S}_D^{1,cr} (\mathcal{T}_h)= \{ v_h \in L^{\infty}(\omega):\ & v_h \vert_T \in \mathcal{P}_1(T) \text{ for all } T \in \mathcal{T}_h \\ &  v_h \text{ continous in } x_S \text{ for all } S \in \mathcal{S}_h  \\ & \text{and } v_h(x_S) = 0 \text{ for all } S \in \mathcal{S}_h \cap \Gamma_D \}
\end{align*}
with the midpoint $x_S$ of side $S \in \mathcal{S}_h$. This means we get the discrete system,
where we search for $v_h \in \mathcal{S}^{1,cr}_{\text{glide}}(\mathcal{T}_h)^2$ and
$p_h \in \mathcal{L}^0(\mathcal{T}_h)$, s.t.
\begin{align*}
\int_\omega v_h \cdot w_h \, \text{d}(r,z) + & \int_\omega \nabla_{\scriptscriptstyle{\mathcal{T}}} v_h : \nabla_{\scriptscriptstyle{\mathcal{T}}} w_h \, \text{d}(r,z) \\- &\int_\omega p_h \text{ div}_{\scriptscriptstyle{\mathcal{T}}}(r w_h) \, \text{d}(r,z) = -J'(\omega, w_h) \\ &\int_\omega q_h \text{ div}_{\scriptscriptstyle{\mathcal{T}}}(rv_h) \, \text{d}(r,z) = 0 
\end{align*}
for all $w_h \in \mathcal{S}^{1,cr}_{\text{glide}}(\mathcal{T}_h)^2$ and $q_h \in \mathcal{L}^0(\mathcal{T}_h)$.

In practice we approximate the weighted integral with a midpoint scheme. This allows us to use the general theory for the Fortin interpolant associated with the Stokes system, which guarantees the well-posedness and stability of the discrete scheme (c.f. \cite{boffi}).

 The discretization of the  Stokes system together with the convexity constraint leads to a minimization problem of the following form
\begin{equation}
\min_{y \in \mathbb{R}^n} 1/2 y^\top Ay -f^\top y \text{ s.t. } By = g \text{ and } Cy \le c
\end{equation}
for  suitable $A \in \mathbb{R}^{n\times n}, B \in \mathbb{R}^{m_1\times n}, C \in \mathbb{R}^{m_2\times n}$ and $f\in \mathbb{R}^n, g \in \mathbb{R}^{m_1}, c \in \mathbb{R}^{m_2}$. This can be formulated as a saddle point problem with an inequality constraint,
\begin{equation}
\min_{y \in \mathbb{R}^n} \max_{z_1 \in \mathbb{R}^{m_1}} \max_{z_2 \in \mathbb{R}^{m_2}, z_2 \ge 0} 1/2 y^\top Ay -f^\top y +z_1^\top(By-g) +z_2^\top(Cy-c).
\end{equation}
This is implemented by including the inequality constraint via a Lagrange multiplier into an Uzawa algorithm, c.f. \cite[Chapter 2.4.3]{glowinski} and Algorithm \ref{algorithm_uzawa}. How to select a suitable stepsize $\alpha$ and a termination criterion as well as extensions to conjugate gradients or with a preconditioner can then be achieved similar to the Uzawa algortihm, c.f. \cite[Section 6.1.5]{pde} or \cite[Section IV.5]{braess}.

\begin{algorithm} 
	\caption{Uzawa algorithm with an inequality constraint} \label{algorithm_uzawa}
	Data:  matrices $A,B,C$ and vectors $f,g,c$ given by the minimization problem\\
	Parameters: stepsize $\alpha$ \\
	Result: minimizer $u$
	\begin{algorithmic}[1]
		\State Set $z^1_0 = 0 \in \mathbb{R}^{m_1}$, $z^2_0 = 0 \in \mathbb{R}^{m_2}$ 
		\For {$k=1,2,\ldots$}
			\State $Au_{k} = f-B^\top z^1_{k-1} -C^\top z^2_{k-1}$
			\State $z^1_k = z^1_{k-1} +\alpha(Bu_k-g)$
			\State $z^2_k = [z^2_{k-1} +\alpha(Cu_k-c)]^+$		
		\EndFor
		
	\end{algorithmic} 
\end{algorithm}

We briefly note that the approach taken in \cite{BW18}, where no constraint on the volume was posed, and the deformation field was computed from
a problem of linear elasticity, did not work well in the problems under consideration, since it resulted in a
poor approximation near the axis, due to the weight $r$ from the transformation. Because of the preservation of volume, this effect occurred only moderately when using the Stokes equation to compute the deformation
field.

\subsection{Line Search}\label{sec:line_search}
We now list the conditions imposed for the deformation field  to find the step
size $\tau>0$ used to update the domain.  We search for the smallest non-negative integer
$k$ such that for $\tau = \tau_0^k$ the following four conditions hold:
\begin{enumerate}
	\item \label{cond:1} The boundary $\Gamma_{\textup{out}}$ avoids
		self-penetration, i.e. the convex curve describing the boundary is
		injective.
\item The linearized convexity constraint is met.
\item  \label{cond:3} The objective value does decrease. 
\item \label{cond:4} The preservation of volume is met, up to a prior set
	tolerance. This was in part necessary, since otherwise the volume was observed to
	change drastically, which makes it difficult to find a suitable stopping
	criterion and to interpret the results. With this condition, the algorithm
	showed better results, but needed more iterations in most cases.
\end{enumerate}
The objective value mentioned in condition \ref{cond:3} is evaluated on the newly generated triangulation, rather than the deformed triangulation. Formally, this means that the line search might not terminate.
 However in practice, this way the shape optimization algorithm needed less iterations to find a stationary domain, since the potential increase of the objective value of the updated domain due to the remeshing of the domain was avoided. No significant difference was observed for the approximated optimal domains and optimal values, if the line search was performed on the deformed triangulation.
 
The algorithm terminates if either $\vert J'(\omega_h,V_h)\vert <
\varepsilon_{\text{stop}}$ or if $\tau < \tau_{min}$. In practice, the latter was usually the reason for termination, due to the second and third condition of the line search, i.e. the objective value did no longer decrease under the convexity constraint. In general, this was observed for either option for the comparison of the objective value in condition \ref{cond:3}.\color{black}

\subsection{An Implementable Code}
\label{sec:pseudo_code}

The following algorithm \ref{algorithm} illustrates the conceptual design of our code, based on \cite{BW18}.

\begin{algorithm} 
	\caption{Shape Optimization Algorithm} \label{algorithm}
	Data: boundary curve $\Gamma_{\textup{out}}^{h}$ of initial domain, the objective functional $J$ to minimize \\
	Parameters: initial step size $\tau_0 >0$, convergence tolerance
	$\varepsilon_{tol} >0$, minimal step size $\tau_{min}$  \\
	Result: boundary curve  $\Gamma_{\textup{out}}^{h}$ of improved domain 
	\begin{algorithmic}[1]
		\State Generate a triangulation $\mathcal{T}_h$ of $\omega_h$ to $\Gamma_{\textup{out}}^{h}$
		\For {$i=1,2,\ldots$}
			\State Approximate shape gradient $J^\prime(\omega_h,\cdot)$
			\State Calculate deformation field  $V_h$ under linearized convexity constraint
			\If {$\vert J_h^\prime(\omega_h,V_h) \vert \le \varepsilon_{tol}$}
			\State STOP, the current iterate $\omega_h$ is almost stationary;
			\EndIf 
			\State Set $k = 0$;
			\While {Condition \ref{cond:1} to \ref{cond:4} are violated for $\tau = \tau_0^k$}
				\State $k = k+1$
			\EndWhile
			\State $\tau = \tau_0^k	$
			\If {$\tau < \tau_{min}$}
			\State STOP, the line search failed;
			\EndIf 
			\State Move the boundary curve according to $\Gamma_{\textup{out}}^{h} = (I + \tau V_h)(\Gamma_{\textup{out}}^{h})$
			\State Generate triangulation $\mathcal{T}_h$ of $\omega_h$ to updated boundary curve $\Gamma_{\textup{out}}^{h}$
		\EndFor
		
	\end{algorithmic} 
\end{algorithm}
\section{Numerical Experiments}
\label{chap:c_num_ex}

\subsection{First Eigenvalue of the Dirichlet Laplacian} \label{chap:num_dirichlet}

For the first eigenvalue of the Dirichlet Laplacian, it is well known that the
optimal domain among open, convex shapes of a certain volume is the ball, see
\cite{krahn1} and \cite{krahn2}. 
Therefore we will use this example to validate the shape optimization
algorithm, by looking at the results for different initial domains and mesh
sizes.

\begin{figure}[p]
\begin{center}
    \includegraphics[height = 6cm]{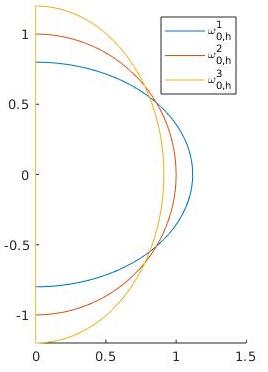} \quad
     \includegraphics[height = 6cm]{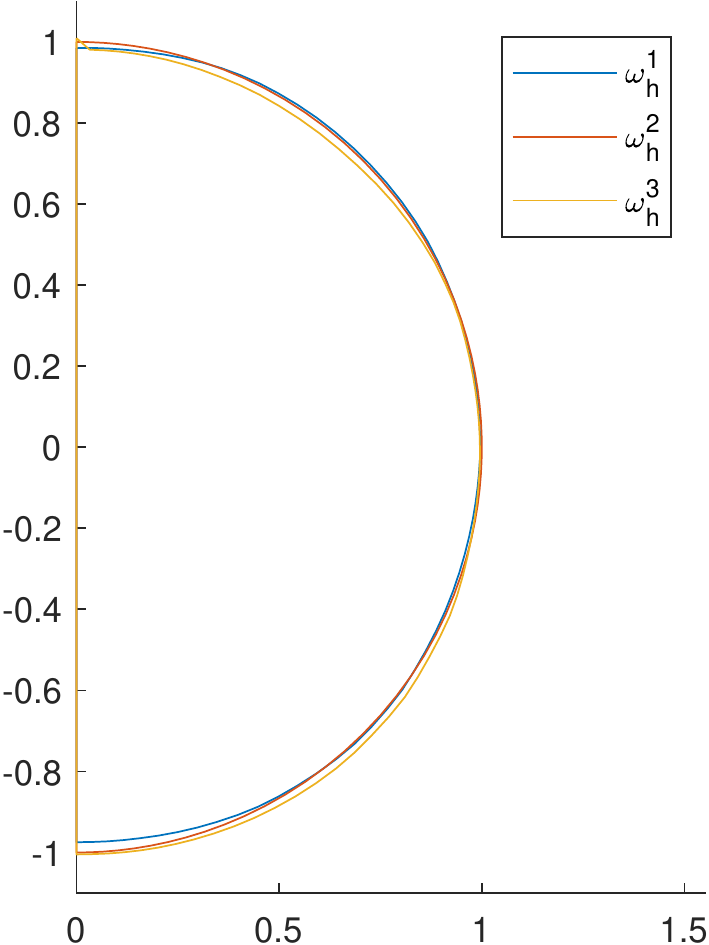}
     \end{center}
  \caption{Initial domains $\omega_{0,h}^i$ and resulting optimal domains $\omega_h^i$, $i = 1,2,3$ of \eqref{problem_dirichlet}, approximately a ball, with constraint $\vert \omega \vert_r = 2/3$}
  \label{fig:dir_tr}
    \end{figure}

Similar to the eigenvalue problem in Section \ref{chap:ex_theory} we can derive the rotationally reduced two-dimensional eigenvalue problem:
\begin{align*}\label{problem_dirichlet}
\textup{Minimize } & \lambda_1(\omega) = J(u,\omega) = \int_{\omega} \vert \nabla u \vert ^2 r \, \text{d}(r,z) \\ 
 \tag{$\mathbf{P_D}$}\textup{w.r.t } & \omega \subset Q \subset \mathbb{R}_+\times \mathbb{R} \text{ open and convex and }2\pi\vert \omega \vert_r = M \\
\text{s.t. } & u \in H^1_r(\omega) \text{ with } \Vert u \Vert_{L^2_r} = 1 \text{ is eigenfunction to the reduced problem} \\ &\begin{cases}
 -(\partial_r u +   r \partial_r^2 u +    r \partial_h^2 u)=  \lambda_1 r u &\text{ in } \omega \\
 u = 0 &\text{ on } \Gamma_{\textup{out}}
\end{cases} 
 \\ \text{ and } & \text{the rotated domain  } R(\omega) \subset  \mathbb{R}^3 \text{ is also convex}.
\end{align*} 
The shape optimization was executed for different mesh refinements and initial
domains. Chosen as initial domains were half-ellipsoids with radii $(a_i,r_i)$
with $a_1 = 0.8$, $a_2 = 1$ and $a_3 = 1.2 $ and $r_i$ so that $\vert
\omega_{0,h}^i \vert_r = 2/3$ for $i = 1,2,3$, so that the volume of the
corresponding three-dimensional domain is the same as that of the unit ball.
The approximated eigenvalues are listed in Tables \ref{table:dir1_h_1_tr},
\ref{table:dir1_h_2_tr} and \ref{table:dir1_h_3_tr} and the initial and approximated optimal
domains for $h = 2^{-5}$ can be seen in Figure \ref{fig:dir_tr}. For reference,
$\lambda_1(B_1(0)) = j^2_{3/2-1,1}  = \pi^2 \approx 9.8696 $, c.f. \cite[(1.13)]{henrot}, and the approximated eigenvalue
$\lambda_1^h(\mathcal{T}_h(B_1(0)) \approx 9.8753 $ for $h=2^{-5}$. The experimental results show that the optimal value known from the Faber-Krahn inequality is approximated well, and suggest a linear rate of convergence, see Tables \ref{table:dir1_h_1_tr} to \ref{table:dir1_h_3_tr}. The error in the preservation of volume for refinements of $h \le 2^{-3}$ is below $10^{-2}$.

\begin{table} [p]
\begin{center}
\begin{tabular}{ l | c | c | c | c| c }
 & $\lambda_1^h(\omega_{0,h}^1)$ &  $\vert\omega_{0,h}^1\vert_r$ & $\lambda_1^h(\omega_{h}^1)$ &  $\vert\omega_{h}^1\vert_r$ & $\vert \lambda_1^h(\omega_{h}^1) - \lambda_1^h(\mathcal{T}_h(B_1(0)) \vert $\\ \hline
 $h = 2^{-2}$ & 10.5403  &  0.6381   
  & 9.9777 & 0.6391 & 0.1081  \\
 $h = 2^{-3}$ & 10.4308  &  0.6593   
  &  9.9624 & 0.6572 &  0.0928 \\
$h = 2^{-4}$ & 10.3785   & 0.6648
  & 9.9231  & 0.6617 &   0.0535\\
$h = 2^{-5}$ & 10.3634  &  0.6662
  & 9.9052 & 0.6634  &     0.0356\\
\end{tabular} 

\caption{Discrete eigenvalues for initial domain $\omega_{0,h}^1$ and resulting optimal domain $\omega_{h}^1$ of \eqref{problem_dirichlet} for different levels of refinement, with constraint $\vert \omega \vert _r=  2/3$  and absolute errors}
\label{table:dir1_h_1_tr}
\end{center}
\begin{center}
\begin{tabular}{ l | c | c | c | c| c  }
 & $\lambda_1^h(\omega_{0,h}^2)$ &  $\vert\omega_{0,h}^2\vert_r$ & $\lambda_1^h(\omega_{h}^2)$ &  $\vert\omega_{h}^2\vert_r$ &$\vert \lambda_1^h(\omega_{h}^2) - \lambda_1^h(\mathcal{T}_h(B_1(0)) \vert $  \\ \hline
 $h = 2^{-2}$ & 10.0218  &  0.6381
  & 9.9054 &  0.6511 & 0.0358 \\
 $h = 2^{-3}$ &  9.9422  &  0.6593 
  & 9.9358 &  0.6601 &  0.0662 \\
$h = 2^{-4}$ & 9.8913    &  0.6648
  & 9.8913  & 0.6648& 0.0217\\
$h = 2^{-5}$ & 9.8753  &  0.6662
  & 9.8753 & 0.6662& 0.0057\\
\end{tabular} 
\caption{Discrete eigenvalues for initial domain $\omega_{0,h}^2$ and resulting optimal domain $\omega_{h}^2$ of \eqref{problem_dirichlet} for different levels of refinement, with constraint $\vert \omega \vert _r=  2/3$  and absolute errors}
\label{table:dir1_h_2_tr}
\end{center}
\begin{center}
\begin{tabular}{ l | c | c | c | c  | c}
 & $\lambda_1^h(\omega_{0,h}^3)$ &  $\vert\omega_{0,h}^3\vert_r$ & $\lambda_1^h(\omega_{h}^3)$ &  $\vert\omega_{h}^3\vert_r$ &$\vert \lambda_1^h(\omega_{h}^3) - \lambda_1^h(\mathcal{T}_h(B_1(0)) \vert $\\ \hline
 $h = 2^{-2}$& 10.3186  &  0.6381 
  & 9.8799  & 0.6575  &  0.0103 \\
 $h = 2^{-3}$ & 10.2321  &  0.6593
  & 9.9332 & 0.6593& 0.0636 \\
$h = 2^{-4}$ &  10.1732  &  0.6648
  & 9.8990 & 0.6642  &  0.0294\\
$h = 2^{-5}$& 10.1578  &  0.6662 
  &  9.8917 & 0.6645 &  0.0221 \\
\end{tabular}
       
\caption{Discrete eigenvalues for initial domain $\omega_{0,h}^3$ and resulting optimal domain $\omega_{h}^3$ of \eqref{problem_dirichlet} for different levels of refinement, with constraint $\vert \omega \vert_r=  2/3$ and absolute errors}
\label{table:dir1_h_3_tr}
\end{center}
\end{table}

\subsection{Eigenvalue Arising in a Problem of Optimal Insulation}
The  reduced variant of the problem of optimal insulation led to the following two-dimensional discrete problem.

\begin{align*}\label{problem_m_discrete}
\textup{Minimize } &  \lambda^{r,\varepsilon,h}_m(\omega_h) \\
\textup{w.r.t. }&\omega_h \subset \mathbb{R}_+\times \mathbb{R}, \ \mathcal{T}_h \in \mathbb{T}_{c_{\textup{usr}}} \textup{ triangulation of } \omega_h \tag{$\mathbf{P_{m,h,\varepsilon}}$} \\
\textup{s.t. }& u_h  \in \mathcal{S}^1(\mathcal{T}_h) \textup{ solves } \eqref{eq:reg_var_dis} \textup{ and } \Vert u_h \Vert_{L^2_r(Q)} = 1
 \\ & \omega_h \subset Q \textup{ is  convex and open and } 2\pi\vert \omega_h \vert_r = M \\ \text{and } & \text{the rotated domain  } R(\omega_h) \subset  \mathbb{R}^3 \text{ is also convex} 
\end{align*}
for the class $\mathbb{T}_{c_{\textup{usr}}}$  of conforming, uniformly
shape regular triangulations $\mathcal{T}_h$ of polyhedral subsets of
$\mathbb{R}^2$ with $h_T/\varrho_T \le c_{\textup{usr}}$  for all
elements $ T \in \mathcal{T}_h$ with diameter $h_T \le h$ and inner radius $\varrho_T$ for a universal
constant $c_{\textup{usr}} >0$. 

We look at several values for the mass $m$. From \cite{BBN17} we know,
that for the ball symmetry breaking for the distribution of insulating material
occurs if $m$ is below a critical value.
\begin{theorem}[c.f. \cite{BBN17}, Theorem 3.1] \label{thm:sym_breaking}
Let $\Omega$ be a ball. Then there exists $m_0>0$ such that the eigenfunction
to \eqref{eq:rayleigh_oi} is radial if $m > m_0$, while the solution is not
radial for $0 < m <m_0$. As a consequence, the optimal insulation thickness
$\ell_{opt}$ is not constant if $m <m_0$.
\end{theorem} 
In \cite{BBN17} it is further noted, that this threshold is given by the
unique positive $m$ for which $\lambda_m = \mu_2$, the first non-zero
eigenvalue of the Neumann problem. Furthermore it is proven, that for $m <
m_0$ the ball is not a stationary domain for the shape optimization problem.
We can use the Neumann eigenvalue to approximate the value for the threshold
$m_0$ of the dimensionally reduced problem for the ball, which is given by
approximately $m_0\approx 5.7963$.

\begin{figure}[h]
\begin{center}
     \includegraphics[height = 4cm]{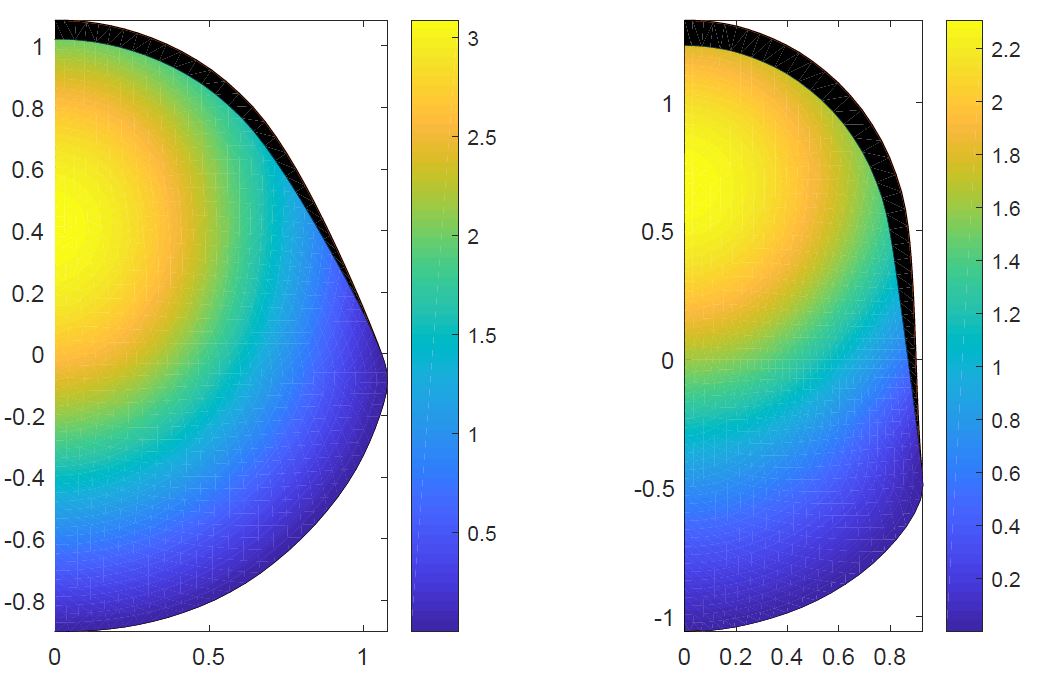}
  \caption{Optimal domains $\omega_h$ of \eqref{problem_m_discrete} for $m =2$ (left), $m=5$ (right), with boundary film (black, scaled with $\epsilon= 1/10$) and eigenfunction (shaded), with constraint $\vert \omega \vert_r = 2/3$; corresponding eigenvalues c.f. Table \ref{table:oi_tr1}}
  \label{fig:oi_tr_2}
   \label{fig:oi_tr_5}
\end{center}
    \end{figure}

\begin{figure}[h] 
 \begin{center}
 \includegraphics[height = 3.5cm]{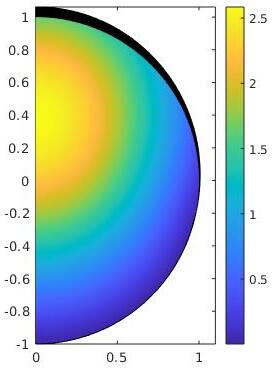} \quad
 \includegraphics[height = 3.5cm]{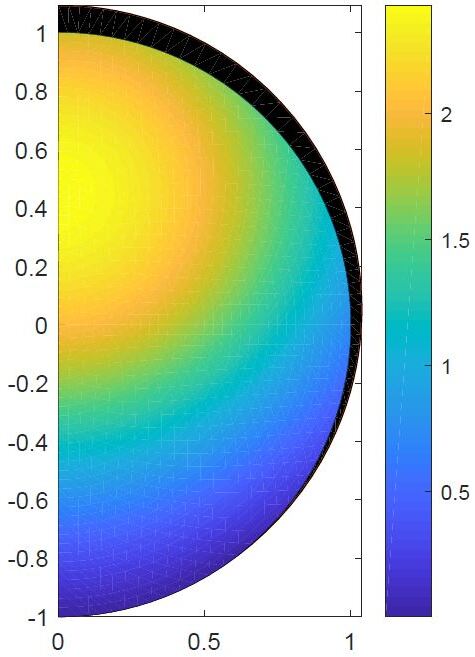}
  \caption{Non-radial eigenfunctions (shaded) with boundary film (black, scaled with $\epsilon= 1/10$) for $m = 2$ (left) and $m = 5$ (right) for the fixed ball}
  \label{fig:non_radial_ball}
 \end{center}
    \end{figure}

\begin{figure}[h]
\begin{center}
\includegraphics[width = 2.5cm]{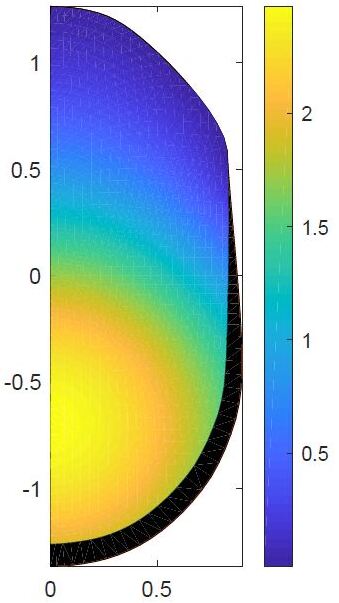} \quad 
\includegraphics[width = 2.5cm]{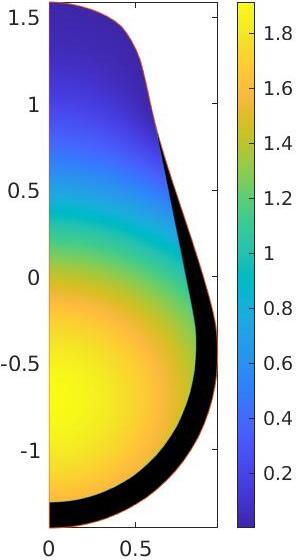} \quad 
\includegraphics[width = 2.5cm]{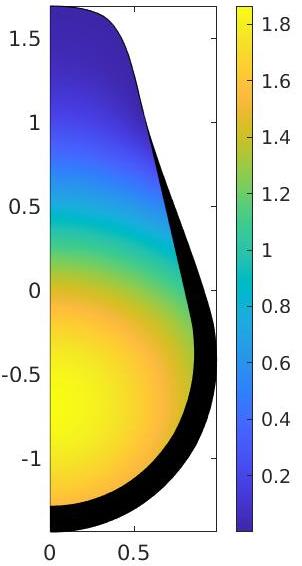} \quad
\includegraphics[width = 2.5cm]{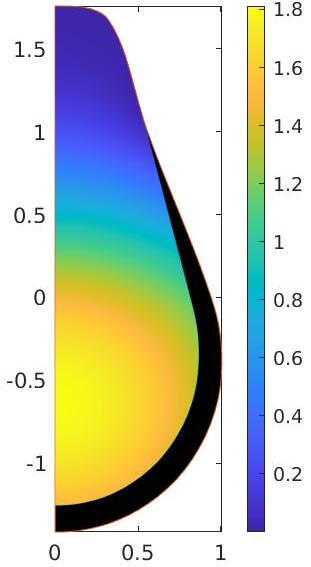}\\
 \includegraphics[width = 2.5cm]{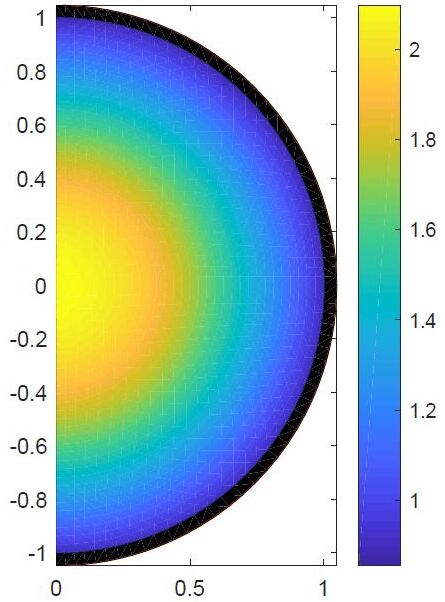}  \quad  \includegraphics[width = 2.5cm]{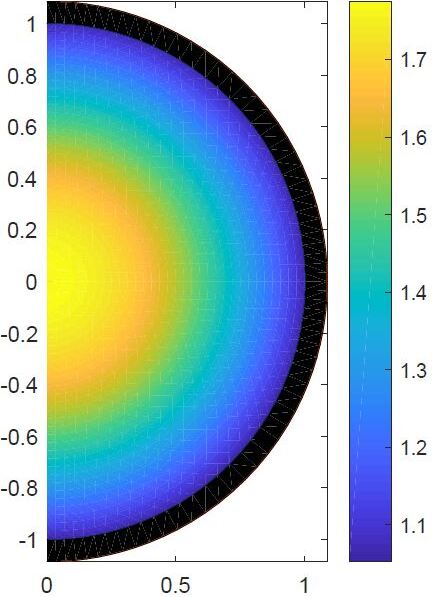} \quad  \includegraphics[width = 2.5cm]{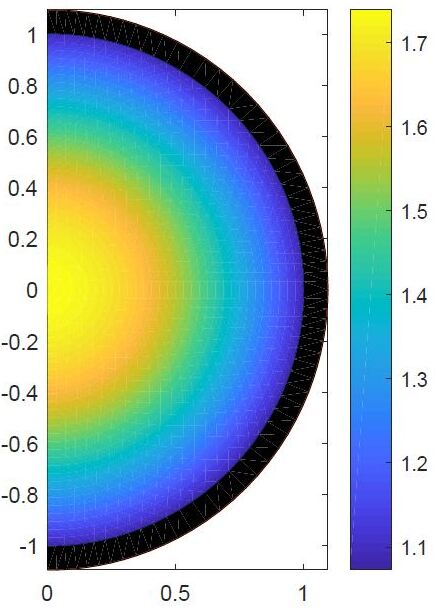}  \quad 
 \includegraphics[width = 2.5cm]{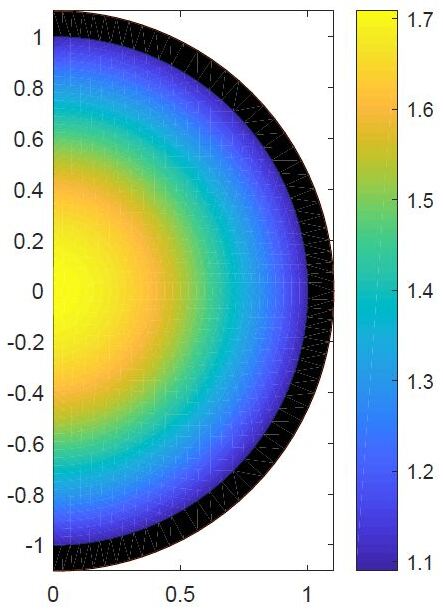}
  \caption{ Approximated stationary domains $\omega^i_h$ ($i = 1$ top, $i = 2$ bottom) of \eqref{problem_m_discrete} for $m =6,11,12$ and $13$ (left to right), with boundary film (black, scaled with $\epsilon= 1/10$) and eigenfunction (shaded), with constraint $\vert \omega \vert_r = 2/3$; corresponding eigenvalues c.f. Table \ref{table:oi_tr2}\color{black}. For $m = 6,11$ the asymmetric domains (top) are optimal, for $m = 12,13$ the half-disk (bottom) is optimal. The other domains are stationary but not optimal under the shape optimization.}
  \label{fig:oi_tr_6}
\end{center}
    \end{figure}

We next address numerical approximations for the values $m = 2,5,6,11,12$ and $13$.  The experimental results displayed in Tables \ref{table:oi_tr1} and \ref{table:oi_tr2} and Figures \ref{fig:oi_tr_2}, \ref{fig:non_radial_ball}  and  \ref{fig:oi_tr_6} were obtained on triangulations  $\mathcal{T}_h$ with maximal mesh size $h=2^{-5}$ and regularization parameter  $\varepsilon = N^{-1/2}/10$, where $N$ is the number of nodes of $\mathcal{T}_h$. The numerical experiments confirm that for the two values lower then the critical mass, the ball is no stationary domain. Only one asymmetric optimal domain was found for each value of $m$, c.f. Figure \ref{fig:oi_tr_2}
and Table \ref{table:oi_tr1}. For the larger values, in the numerical experiments the ball is also
stationary, and for $m = 11,12$ and $13$ it is experimentally optimal, see Figures \ref{fig:oi_tr_6} and Tables \ref{table:oi_tr2}. As in Section \ref{chap:num_dirichlet} half-ellipsoids with different ratios were chosen as initial domains. For values of $m$ where more than one stationary domain was approximated, the result of the optimization algorithm depended on the choice of the ratios for the initial domains. In general, depending on the value $m$, when initial domains were chosen that are more prolate, an asymmetric domain was approximated, while oblate ellipsoids and ellipsoids closer to the ball led to the ball being approximated.

The algorithm only detects local minima with the approximated domains depending on the initial domains, so the stationary domains approximated might not be global solutions. \\ 
When comparing the approximated optimal domains with each other, we also observe that for the non-radial solutions, a large portion of the
insulating film concentrates in one area, which creates a hotspot inside the domain,
where the temperature is preserved better, while other areas are neglected,
having no insulating material on the boundary. 

\begin{remark} \rm
We briefly note, that even though we search for eigenfunctions among
rotationally symmetric functions, the numerical results are still consistent
with the expectations we have from \cite{BBN17} regarding the radial
symmetry for the eigenfunctions for the ball. We were able to observe that for
$m < m_0$, the critical value related to the Neumann eigenvalue, c.f. Theorem
\ref{thm:sym_breaking}, the eigenfunction is no longer positive or symmetric,
c.f. Figure \ref{fig:non_radial_ball}, but for values $m > m_0$ it is. 
Further, the shape optimization problem showed, that for the chosen values $m <
m_0$ the ball was no stationary domain, while for $ m>m_0 $ it was stationary
and for the higher values even optimal. This consistency of the numerical
results suggests, that the restriction to rotationally symmetric functions is
justified.  \\ The critical value $m_0$ relates to the symmetry of the eigenfunctions on the ball and whether the ball is a stationary domain. 
Theorem 3.1 in \cite{BBN17} does not consider the optimality of the ball under the shape optimization. However, the experimental results of the shape optimization suggest that there might be another critical value of mass $m_1$,  such that for $m < m_1$ an asymmetric domain yields an optimal eigenvalue, while for $m > m_1$ the ball is the optimal domain.
\end{remark} 

We want to take a closer look at the properties of the optimal domains for the
eigenvalue problem in optimal insulation. First we will look at the improvement
of the eigenvalue the shape optimization provides and afterward at the optimal
domains themselves. The experiments in this section were obtained with a triangulation with a maximal mesh size $h = 2^{-4}$ and $\varepsilon$ chosen as in the previous experiments.

\begin{figure}[h]
\begin{center}
\includegraphics[height = 5cm]{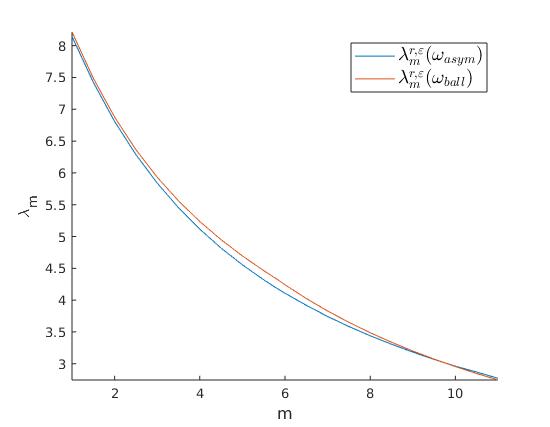}\includegraphics[height = 5cm]{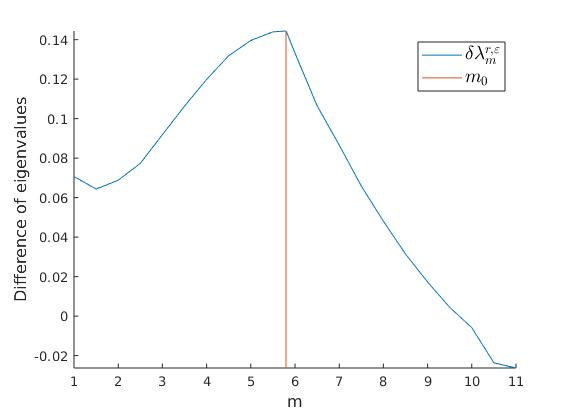}
\caption{Eigenvalues of the asymmetric stationary domain and the ball compared for different values of $m$ (left) and difference in eigenvalues $\delta \lambda_m^{r,\varepsilon} = \lambda_m^{r,\varepsilon}(\omega_{\text{ball}})-\lambda_m^{r,\varepsilon}(\omega_{\text{asym}})$
with peak at the critical value $m_0$ (right)}
\label{fig:comp_ev}
\end{center}
\end{figure}

Comparing the eigenvalue of the ball to that of the respective stationary
asymmetric domain for different values of mass $m$, c.f. Figure
\ref{fig:comp_ev}, shows that the benefit of the shape optimization is
greatest around the critical value $m_0$.

\begin{figure}
\begin{center}
\includegraphics[height = 4cm]{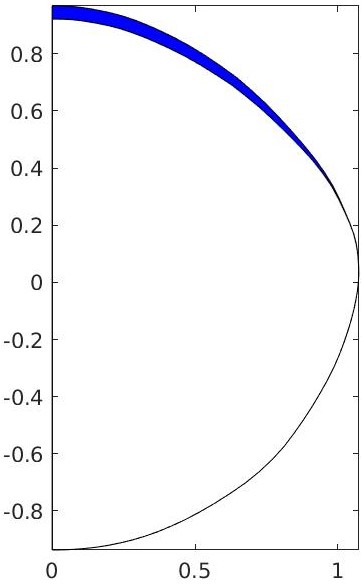}
\includegraphics[height = 4cm]{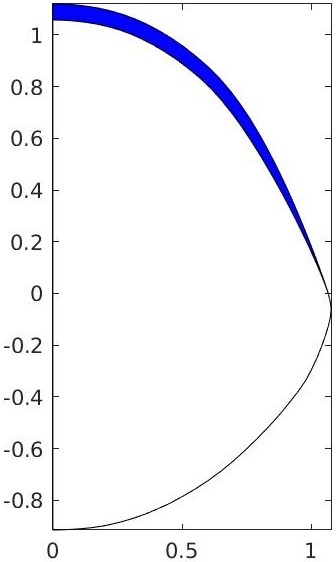}
\includegraphics[height = 4cm]{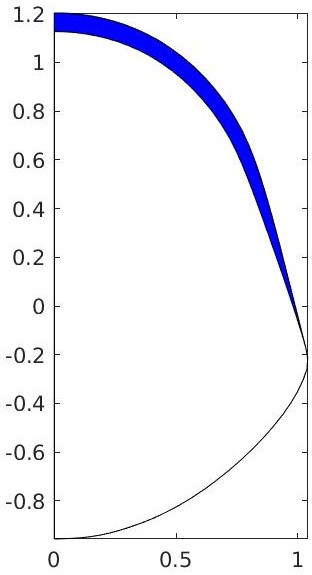}
\includegraphics[height = 4cm]{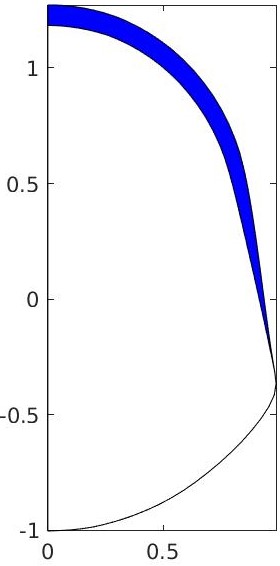} 
\includegraphics[height = 4cm]{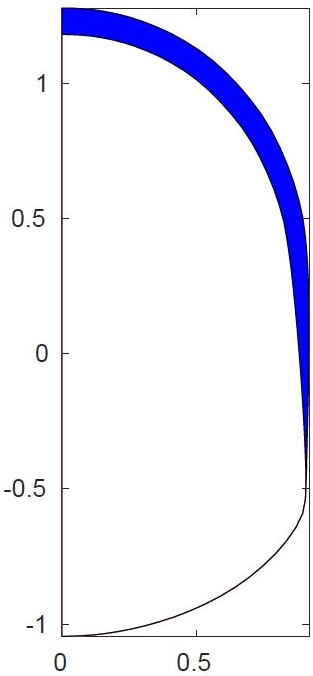} 
\includegraphics[height = 4cm]{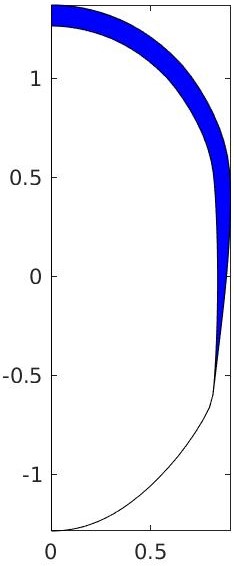}\includegraphics[height = 4cm]{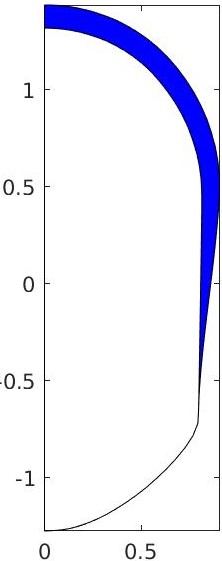}\includegraphics[height = 4cm]{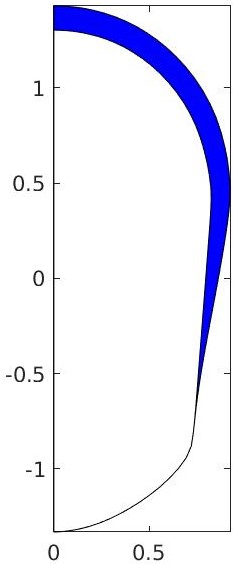} 
\includegraphics[height = 4cm]{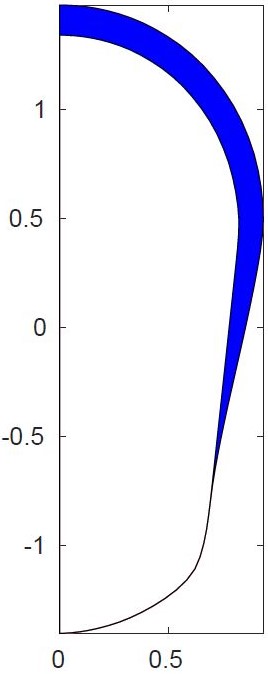}\includegraphics[height = 4cm]{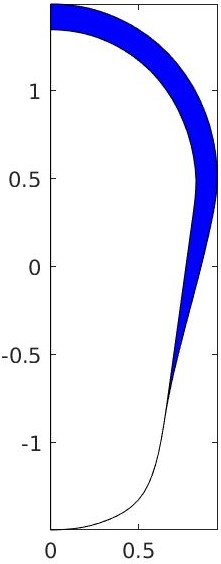}\includegraphics[height = 4cm]{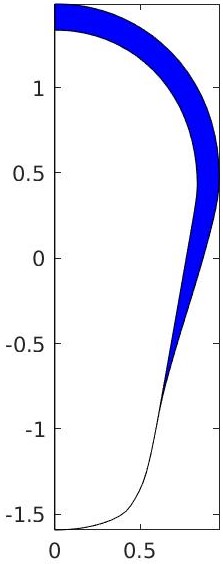}\includegraphics[height = 4cm]{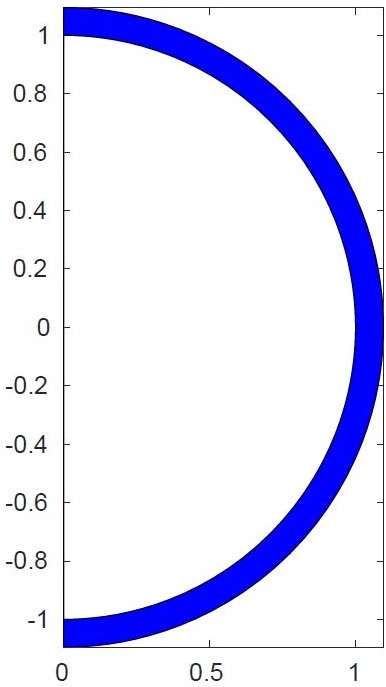}
\caption{Approximated stationary and optimal domains for $\lambda^{r,\varepsilon}_m$ for $m = 1$ to $12$ (left to right, top to bottom) and insulating film (blue, scaled with $\epsilon= 1/10$)}
\label{fig:evol_optimal_domain}
\end{center}
\end{figure}

Next, we take a closer look at the optimal domains. For the masses $ m =
1,2,..,12$ the optimal domains with insulating film are displayed in Figure
\ref{fig:evol_optimal_domain}. As $m$ decreases, the eigenvalue is closer to
the eigenvalue of the Dirichlet Laplacian.  We notice, for $m = 1$, the approximated optimal
shape is closer to a ball, and as the values $m$ increase the optimal domains
become more prolate, until, for $m = 12$, the ball is the approximated optimal domain. 

\begin{figure}
\begin{center}
\includegraphics[height = 4cm]{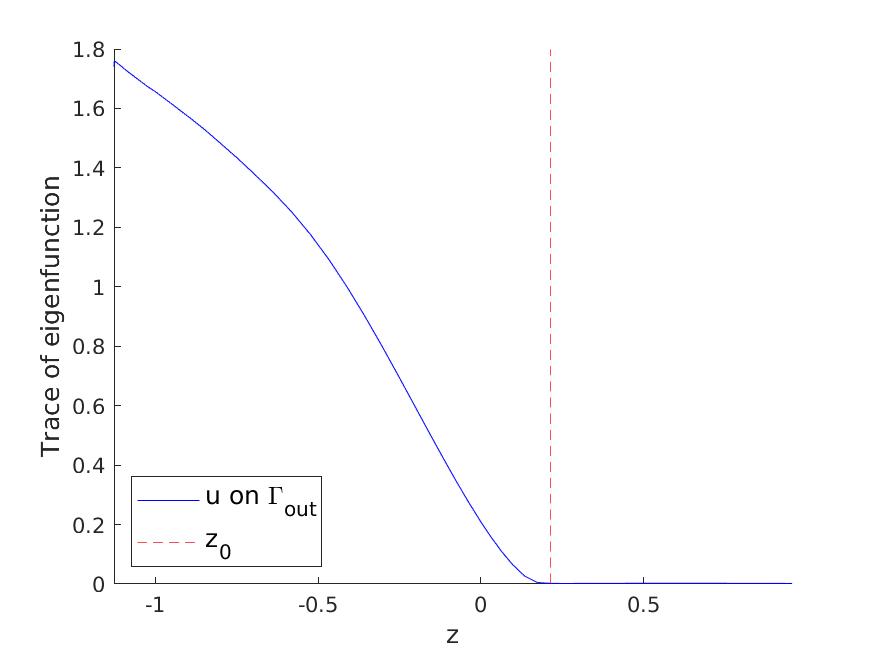} \hspace{0.5cm} 
\includegraphics[height = 4cm]{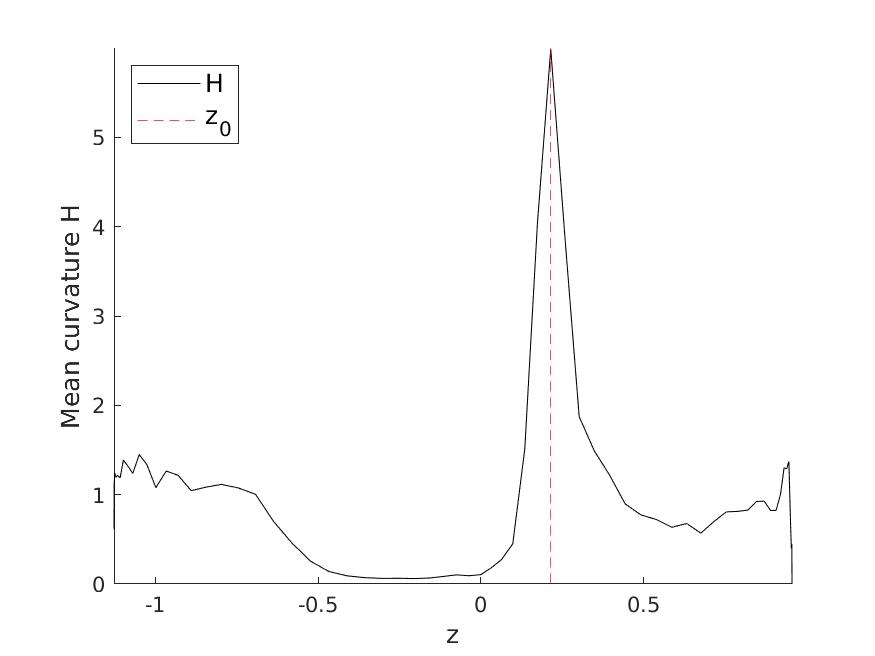} \\
\includegraphics[height = 4cm]{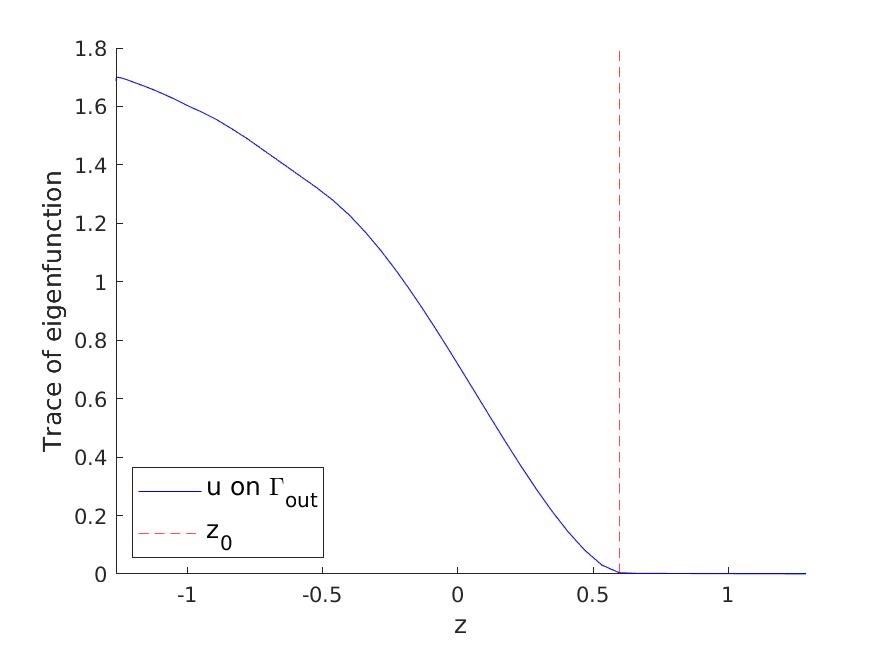}\hspace{0.5cm} 
 \includegraphics[height = 4cm]{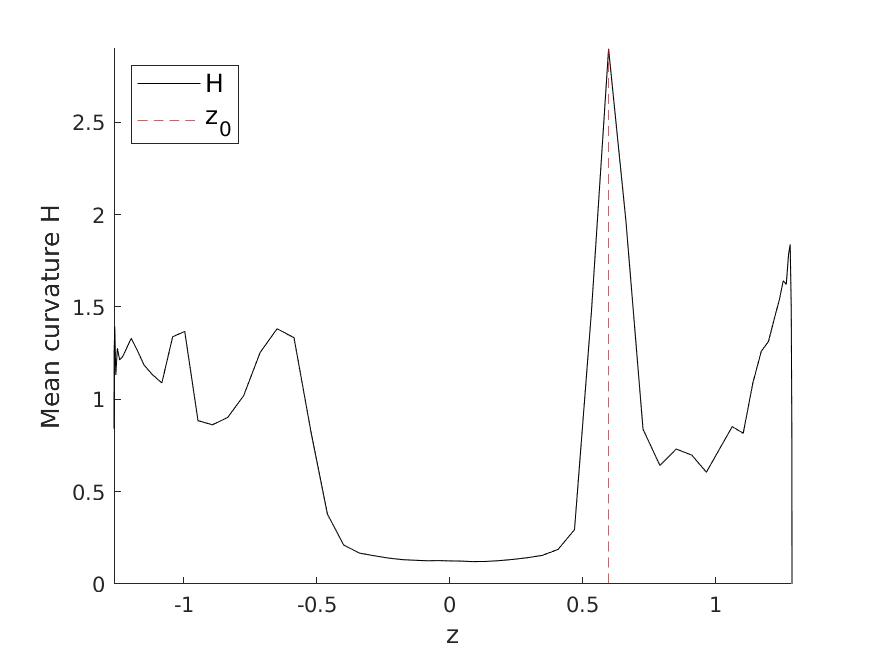} \\
\includegraphics[height = 4cm]{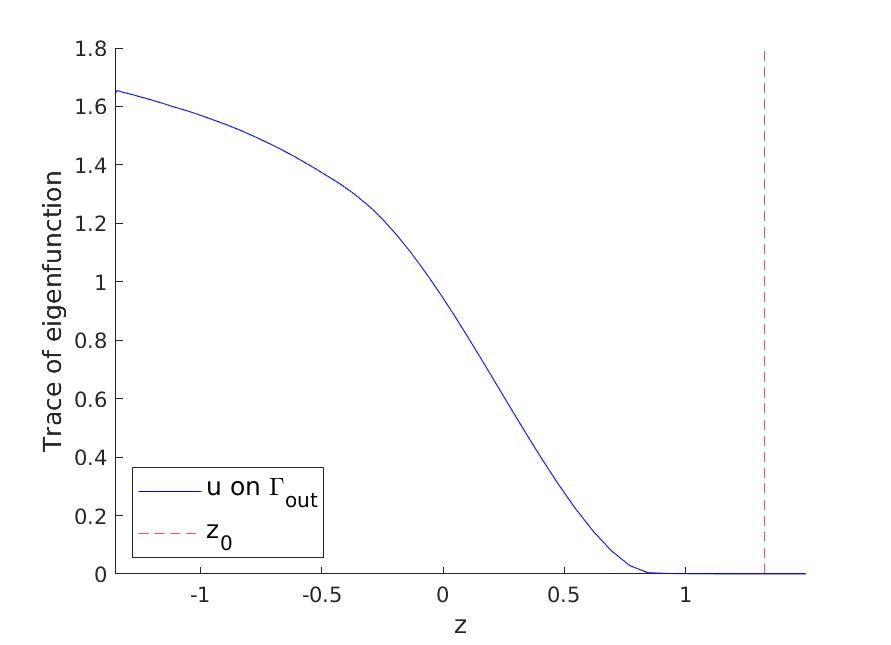} \hspace{0.5cm} \includegraphics[height = 4cm]{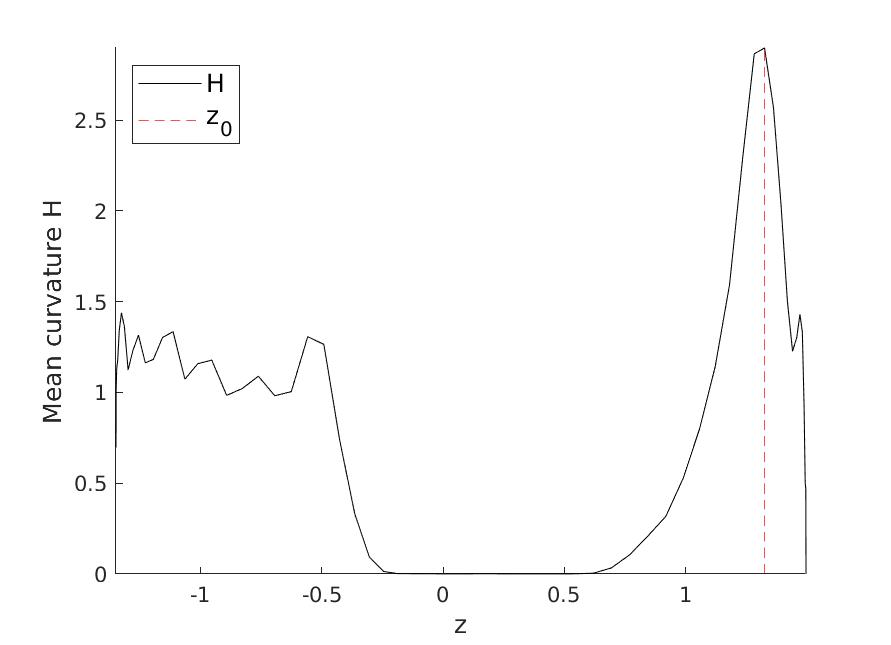}
\caption{Values of eigenfunctions (left) and mean curvature (right) along the boundary curve (as a function of $z$ along the axis of rotation) of optimal domains for $m = 3,6$, and $10$ (top to bottom), with approximate location $z_0$ of the kink in the boundary on a triangulation with maximal mesh size $h = 2^{-4}$. The corresponding domains are shown in Figure \ref{fig:evol_optimal_domain}. }
\label{fig:u_bdy}
\end{center}
\end{figure}

We further notice that for the asymmetric domains a kink is
formed around the surface area where the eigenfunction is zero. For $m = 2,...,9$, this kink might even be non-smooth. By taking a closer look at the
values of the eigenfunction at the boundary nodes and the mean curvature of the
boundary $\Gamma_{\textup{out}}$ for the approximated optimal domains, c.f. Figure
\ref{fig:u_bdy}, we can see that this kink is located around the surface, where
no insulating material is placed, and that for the values where the kink might
be non-smooth, it is located where the eigenfunction vanishes. The
insulating material then focuses on one side of the kink. The corresponding
domains are those shown in Figure \ref{fig:evol_optimal_domain}. In summary, the numerical experiments suggest, that
the asymmetric optimal domains tend to be non-smooth for lower values of $m$. \\

\textit{Acknowledgments:} This work is supported by DFG grants BA2268/4-2 within the Priority Program SPP 1962 (Non-
smooth and Complementarity-based Distributed Parameter Systems: Simulation and
Hierarchical Optimization).

\begin{table}[h]
\begin{center}
\begin{tabular}{ l | c  | c }
 & $\lambda^{r,\varepsilon}(\omega_{h})$ & $\vert \omega_h \vert_r$\\ \hline
 $m = 2$ & 6.819940118007397 &   0.6595 \\ 
$m = 5$ &4.554732496286795 & 0.6596\\ 
\end{tabular} 
\caption{Eigenvalues $\lambda^{r,\varepsilon}_m$ of approximated optimal domains for different values of $m < m_0$ c.f. Figure \ref{fig:oi_tr_2}, with constraint $\vert \omega \vert_r = 2/3$}\label{table:oi_tr1}
\end{center}
\end{table}

    \begin{table}[h]
\begin{center}
\begin{tabular}{ l | c  | c | c | c  }
 & $\lambda^{r,\varepsilon}_m(\omega^1_{h})$ & $\vert \omega^1_h \vert_r$ & $\lambda^{r,\varepsilon}_m(\omega^2_{h})$ & $\vert \omega^2_h \vert_r$\\ \hline
$m = 6$ &4.112232986601394 &0.6631 & 4.241084607303154 & 0.6662\\
$m = 11$ &2.769366507533780 & 0.6633 &2.744289963928029& 0.6662\\
$m=12$ & 2.606239519805330  & 0.6621& 2.561023079370016 & 0.6662  \\   
$m =13$ & 2.459976343090586 & 0.6633 &2.400341990779929 &0.6662 
\end{tabular} 
\caption{ Eigenvalues $\lambda^{r,\varepsilon}_m$ of approximated stationary domains $\omega^1_h$ (asymetric) and $\omega_h^2$ (half-disk), c.f. Figure \ref{fig:oi_tr_6} for different values of $m$, with constraint $\vert \omega \vert_r = 2/3$}\label{table:oi_tr2}
\end{center}
\end{table}

\end{document}